\documentclass[a4paper]{amsart}                

\usepackage[latin1]{inputenc}                  
\usepackage[T1]{fontenc}                       
\usepackage[spanish,english]{babel}            
\usepackage[pdftex]{graphicx}                  
\usepackage{amsmath,amssymb,amsthm}            
\usepackage{latexsym}                          
\usepackage{delarray}                          
\usepackage{bbm}                               
\usepackage{hyperref}                          
\usepackage[pdftex,usenames,dvipsnames]{color} 
\usepackage{calrsfs,eufrak,mathrsfs,upgreek}   
\usepackage{bbding,trfsigns}                   
\usepackage{wasysym}                           
\usepackage{datetime}                          

\DeclareMathAlphabet{\mathpzc}{OT1}{pzc}{m}{it} 

\newtheorem{Th}{Theorem}[section]              

\newtheorem{Rm}{Remark}[section]

\newtheorem{Prop}{Proposition}[section]

\title[Variation operators on $BMO$ in the Schr\"odinger setting]{Variation operators for semigroups and Riesz transforms on $BMO$ in the Schr\"odinger setting}
\author{J.J. Betancor}
\address{Departamento de Análisis Matemático\\
Universidad de la Laguna\\
Campus de Anchieta, Avda. Astrofísico Francisco Sánchez, s/n\\
38271 La Laguna (Sta. Cruz de Tenerife), Spain}
\email{jbetanco@ull.es}

\author{J.C. Fariña}
\address{Departamento de Análisis Matemático\\
Universidad de la Laguna\\
Campus de Anchieta, Avda. Astrofísico Francisco Sánchez, s/n\\
38271 La Laguna (Sta. Cruz de Tenerife), Spain}
\email{jcfarina@ull.es}

\author{E. Harboure}
\address{Instituto de Matemática Aplicada del Litoral, IMAL\\
Universidad Nacional del Litoral\\
C/ G\"emes\\
Santa Fe, Argentina}
\email{}

\author{L. Rodr\'{\i}guez-Mesa}
\address{Departamento de Análisis Matemático\\
Universidad de la Laguna\\
Campus de Anchieta, Avda. Astrofísico Francisco Sánchez, s/n\\
38271 La Laguna (Sta. Cruz de Tenerife), Spain}
\email{jrguez@ull.es}

\thanks{This paper is partially supported by MTM2007/65609.}

\begin{document}

  \maketitle                                  

  \begin{abstract}
    In this paper we prove that the variation operators associated with the heat semigroup and Riesz transforms related to the Schr\"odinger operator are bounded on the suitable $BMO$ type space.
  \end{abstract}

  \section{Introduction}

  We consider the Schr\"odinger operator $\mathcal{L}$ defined by
  $\mathcal{L}=-\Delta+V(x)$ on $\mathbb{R}^n$, $n\ge 3$. Here $V$
  is a nonnegative and not identically zero function satisfying, for some $q\ge n/2$, the
  following reverse H\"older inequality:

  $(RH_q)$ There exists $C>0$ such that, for every ball $B\subset
  \mathbb{R}^n$,
  $$
  \Big(\int_B V(x)^qdx\Big)^{1/q}\le C\int_B V(x)dx.
  $$
  We write $V\in RH_q$ when $V$ verifies the property $(RH_q)$. Note
  that if $V$ is a nonnegative polynomial, then $V\in RH_q$ for
  every $1<q<\infty$. Also, if $V_\alpha(x)=|x|^\alpha$, $x\in
  \mathbb{R}^n$, $V_\alpha$ belongs to $RH_q$ provided that $\alpha
  q>-n$. Hence, $V_\alpha\in RH_{n/2}$ when $\alpha>-2$, and
  $V_\alpha\in RH_n$ if $\alpha>-1$.

  Harmonic analysis operators related to Schr\"odinger operator
  (Riesz transforms, maximal operators associated with heat and
  Poisson semigroups for $\mathcal{L}$, Littlewood-Paley
  $g$-functions, fractional integrals,...) have been extensively
  studied in last years. The papers of Shen
  \cite{Sh} and Zhong \cite{Zh} can be considered as the starting points. In \cite{Sh}
  and \cite{Zh}  Riesz transforms in the Schr\"odinger setting were
  studied in $L^p$-spaces. Other operators related to $\mathcal{L}$
  were investigated in $L^p$-spaces in \cite{BHS3}, \cite{GLP},
  \cite{LD1}, \cite{LD2} and \cite{S}, amongst others.

  Dziubanski and Zienkiewicz defined Hardy spaces associated with
  $\mathcal{L}$ (see \cite{DZ1}, \cite{DZ2}, and \cite{DZ3}). A
  function $f\in L^1(\mathbb{R}^n)$ is in
  $H_1^\mathcal{L}(\mathbb{R}^n)$ if and only if
  $W^\mathcal{L}_*(f)\in L^1(\mathbb{R}^n)$, where
  $$
W^\mathcal{L}_*(f)=\sup_{t>0}|W_t^\mathcal{L}(f)|,
$$
and $W^\mathcal{L}=\{W_t^\mathcal{L}\}_{t>0}$ denotes the heat
semigroup generated by $-\mathcal{L}$.

The dual space of $H_1^\mathcal{L}(\mathbb{R}^n)$ was investigated
in \cite{DGMTZ}. This dual space was characterized by the space
$BMO^\mathcal{L}(\mathbb{R}^n)$ of bounded  mean oscillation
functions. A function $f\in L^1_{\rm{loc}}(\mathbb{R}^n)$ is in
$BMO^\mathcal{L}(\mathbb{R}^n)$ provided that there exists $C>0$
such that the following two properties are satisfied:

(i) For every $x\in \mathbb{R}^n$ and $r>0$,
$$
\frac{1}{|B(x,r)|}\int_{B(x,r)}|f(y)-f_{B(x,r)}|dy\le C,
$$
where, as usual, $f_{B(x,r)}=\frac{1}{|B(x,r)|}\int_{B(x,r)}f(y)dy$,
and $|B(x,r)|$ denotes the Lebesgue measure of $B(x,r)$; and

(ii) For every $x\in \mathbb{R}^n$ and $r\ge \gamma(x)$,
$$
\frac{1}{|B(x,r)|}\int_{B(x,r)}|f(y)|dy\le C.
$$
Here, the critical radii $\gamma(x)$ is defined by
$$
\gamma(x)=\sup\Big\{r>0:r^{2-n}\int_{B(x,r)}V(y)dy\le 1\Big\}.
$$
Since $V$ is not identically zero, $0<\gamma(x)<\infty$. The norm $\|f\|_ {BMO^\mathcal{L}(\mathbb{R}^n)}$ of $f\in BMO^\mathcal{L}(\mathbb{R}^n)$ is defined by
$$
\|f\|_ {BMO^\mathcal{L}(\mathbb{R}^n)}=\inf\{C>0:\, \rm{(i)}
\,\,\,\rm{and}\,\,\, \rm{(ii)}\,\,\, hold\}.
$$

In \cite{DGMTZ} the behavior on $BMO^\mathcal{L}(\mathbb{R}^n)$ of
certain maximal operators, Littlewood-Paley $g$ functions and
fractional integrals were studied. Also, the
$BMO^\mathcal{L}(\mathbb{R}^n)$-boundedness properties of the Riesz
transforms have been analyzed in \cite{BHS2}, \cite{DL},
\cite{YYZ} and \cite{Y2}.

Suppose that $\{T_t\}_{t>0}$ is a family of operators defined in
$L^p(\mathbb{R}^n)$, $1\le p<\infty$. If $\rho>2$, the
$\rho$-variation operator associated with $\{T_t\}_{t>0}$,
$V_\rho(T_t)$, is defined by
$$
V_\rho(T_t)(f)(x)=\sup_{\{t_j\}_{j\in \mathbb{N}}\downarrow
0}\Big(\sum_{j=1}^\infty
|T_{t_j}(f)(x)-T_{t_{j+1}}(f)(x)|^\rho\Big)^{1/\rho},
$$
where the supremum is taken over all the real decreasing sequences
$\{t_j\}_{j=1}^n$ that converge to zero. The operator
$V_\rho(T_t)$ is related to the convergence of $T_t$, as $t\to 0^+$,
and it estimates the fluctuations near the origin of the family
$\{T_t\}_{t>0}$.

We consider the linear space $E_\rho$ that consists of all those real functions $F$ defined on $(0,\infty)$ such that
$$
\|F\|_{E_\rho} = \sup_{\{t_j\}_{j\in \mathbb{N}} \downarrow 0} \left(\sum^\infty_{j=1}|F(t_j)-F(t_{j+1})|^\rho\right)^{\frac{1}{\rho}} < \infty,
$$
where the supremum  is taken over all the real decreasing sequence
$\{t_j\}_{j=1}^\infty$ that converge to zero. $
\|_{E_\rho}$ is a seminorm on $E_\rho$. The variation operator $V_\rho(T_t)$ can be rewritten in the following way
$$
V_\rho(T_t)(f)(x)=\|T_t(f)(x)\|_{E_\rho}.
$$

The variation operator $V_\rho$ was introduced in the ergodic
context by Bourgain \cite{Bou} (see also Jones et al. \cite{JKRW}).
In last years many authors have investigated the variation
operator associated to semigroups of operators and singular
integrals (\cite{CJRW1}, \cite{CJRW2}, \cite{CMMTV}, \cite{CLMT},
\cite{GiTo}, \cite{HMMTV}, \cite{JSW} and \cite{JW}). Recently,
Oberlin, Seeger, Tao, Thiele and Wright \cite{OSTTW} have analyzed
the variation norm related to Carleson Theorem.

In a previous paper \cite{BFHR} the authors studied the
$L^p$-boundedness properties of the variation operators for the heat
semigroup $\{W_t^\mathcal{L}\}_{t>0}$ and the family of truncated
Riesz transforms
$\{R_\ell^{\mathcal{L},\varepsilon}\}_{\varepsilon>0}$,
$\ell=1,\cdots,n$, in the Schr\"odinger context. Here our goal
is to study the behavior of the variation operators
$V_\rho(W_t^\mathcal{L})$ and
$V_\rho(R_\ell^{\mathcal{L},\varepsilon})$ in
$BMO^\mathcal{L}(\mathbb{R}^n)$. Previously, we analyze the
variation operators $V_\rho(W_t)$ and $V_\rho(R_\ell^\varepsilon)$
in $BMO(\mathbb{R}^n)$, where $\{W_t\}_{t>0}$ and
$\{R_\ell^{\varepsilon}\}_{\varepsilon>0}$, $\ell=1,\ldots,n$,
represent the classical heat semigroup and truncated Riesz
transforms, respectively. By $BMO(\mathbb{R}^n)$ we denote the well
known space of bounded mean oscillation functions in $\mathbb{R}^n$.

 This paper is organized as follows. In Section 2 we state our
results. The proof of the
theorems are showed in Sections 3 (classical setting) and 4
(Schr\"odinger context).

Throughout this paper by $C$ we denote a positive constant that can
change from one line to another one. Moreover, if $B(x_0,r_0)$ with $x_0\in \mathbb{R}^n$ and $r_0>0$, we define $B^*=B(x_0,2r_0)$ and $B^{**}=(B^*)^*$.

\section{Main results}

As it is well known the heat semigroup $\{W_t\}_{t>0}$ generated by
$-\Delta$ is defined, for every $f\in L^p(\mathbb{R}^n)$, $1\le
p\le\infty$, by
$$
W_t(f)(x)=(4\pi
t)^{-n/2}\int_{\mathbb{R}^n}e^{-|x-y|^2/4t}f(y)dy,\,\,\,x\in
\mathbb{R}^n\,\,\,and\,\,\,t>0.
$$

The $L^p$-boundedness properties of the variation operator
$V_\rho(W_t)$, $\rho>2$, was studied in \cite[Theorem 3.3]{JR} and
\cite[Theorem 1.1]{CMMTV}.

\begin{Th} \label{VarLp}(\cite[Theorem 3.3]{JR} and \cite[Theorem 1.1]{CMMTV})
If $ \rho >2$, the variation operator $V_\rho(W_t)$ is bounded from
$L^p(\mathbb R^n)$ into itself, for every $1<p<\infty$, and from
$L^1(\mathbb R^n)$ into $L^{1,\infty}(\mathbb R^n)$.
\end{Th}

In \cite{CMMTV} it was shown that the variation operator
$V_\rho(W_t)$ is not bounded from $L^\infty(\mathbb{R}^n)$ into
itself. Also, they defined in \cite[Section 5]{CMMTV} a function
$f\in L^\infty(\mathbb{R})$ such that $V_\rho(W_t)(f)(x)=\infty$,
a.e. $x\in \mathbb{R}$, for every $\rho>2$. As it is well known
$L^\infty(\mathbb{R}^n)$ is a subset of the space $BMO(\mathbb{R}^n)$
of bounded mean oscillation functions. The behavior of $V_\rho(W_t)$
in the space $BMO(\mathbb{R}^n)$ is stated in the following.

\begin{Th} \label{Th1BMO}
Let $\rho > 2$. If $f \in BMO(\mathbb R^n)$ and $V_\rho(W_t)
(f)(x) < \infty$, a.e. $x \in \mathbb R^n$, then $V_\rho(W_t)f\in BMO(\mathbb{R}^n)$ and
$\|V_\rho(W_t)f\|_{BMO(\mathbb{R}^n)}\le
C\|f\|_{BMO(\mathbb{R}^n)}.$
\end{Th}

For every $\ell=1,\cdots,n$, the Riesz transform $R_\ell(f)$ of
$f\in L^1(\mathbb{R}^n)$, $1\le p<\infty$, is given by
$$
R_\ell(f)(x)=\lim_{\varepsilon\to
0^+}c_n\int_{|x-y|>\varepsilon}\frac{x_\ell-y_\ell}{|x-y|^{n+1}}f(y)dy,\,\,\,a.e.\,\,\,x\in
\mathbb{R}^n,
$$
where $c_n=\Gamma((n+1)/2)/\pi^{(n+1)/2}$.

The  variation operators for $R_\ell,\,\, \ell=1,\ldots,n$, were
investigated in \cite{CJRW1} and \cite{CJRW2}, where the following
$L^p$-boundedness properties were established.

\begin{Th} \label{Th2Lp}(\cite[Theorem 1.2]{CJRW1} and \cite[Theorem A and Corollary 1.4]{CJRW2}). Let $\ell=1,\ldots,n$.
If $\rho >2$, the variation operator $V_\rho(R_\ell^{\varepsilon})$
is bounded from $L^p(\mathbb R^n)$ into itself, for every
$1<p<\infty$, and from $L^1(\mathbb R^n)$ into $L^{1,\infty}(\mathbb
R^n)$.
\end{Th}

By using transference methods Gillespie and Torrea (\cite[Theorem
B]{GiTo}), have proved dimension free $L^p(\mathbb R^n,|x|^\alpha
dx)$ norm inequalities, for every $1<p< \infty$ and $-1 < \alpha <
p-1$, for variation operators of the Riesz transform
$R_\ell,\;\ell=1,\ldots,n$. The idea developed in the proof of
\cite[Lemma 1.4]{GiTo} allows us to analyze the behavior of the
operators $V_\rho(R_\ell^\varepsilon)$ on the space $BMO(\mathbb
R^n)$. Note that if $\ell=1,\ldots,n$, $f\in BMO(\mathbb{R}^n)$ and
$\varepsilon>0$, then the integral
$\int_{|x-y|>\varepsilon}f(y)\frac{y_\ell-x_\ell}{|y-x|^{n+1}}dy$
may be not convergent when $x\in \mathbb{R}^n$. Indeed,
the function $f(x)=\frac{1}{\log(x+2)}\chi_{(0,\infty)}(x)$, $x\in
\mathbb{R}$, belongs to $L^\infty(\mathbb{R}^n)\subset
BMO(\mathbb{R}^n)$ and the limit
$\lim_{N\to\infty}\int_{\varepsilon<|x-y|<N}\frac{f(y)}{x-y}dy$ does
not exist, for any $x\in \mathbb{R}$ and $\varepsilon>0$. However,
it is clear that, for every $0<\varepsilon<\eta$ and $x\in
\mathbb{R}^n$,
$\int_{\varepsilon<|x-y|<\eta}\frac{|f(y)|}{|x-y|^n}dy<\infty$.
Here, the operators $V_\rho(R_\ell^\varepsilon)$ are defined on
$BMO(\mathbb{R}^n)$ in the obvious way by replacing
$R_\ell^{\varepsilon_j}(f)(x)-R_\ell^{\varepsilon_{j+1}}(f)(x)$ by
$c_n\int_{\varepsilon_{j+1}<|x-y|<\varepsilon_j}f(y)\frac{y_\ell-x_\ell}{|y-x|^{n+1}}dy$,
$\ell=1,\cdots,n$ and $j\in \mathbb{N}$, where $f\in
BMO(\mathbb{R}^n)$. In \cite[Theorem B]{CLMT} it was proved that if
$f\in L^\infty(\mathbb{R}^n)$ and $\rho>2$ then,
$V_\rho(H^\varepsilon)(f)(x)=\infty$, a.e. $x\in \mathbb{R}$, or
$V_\rho(H^\varepsilon)(f)(x)<\infty$, a.e. $x\in \mathbb{R}$, where
$H$ denotes the Riesz transform on $\mathbb{R}$, that is, the
Hilbert transform. Moreover, as it can be seen in \cite[Section
1]{CLMT}, if $f(x)=sgn(x)$, $x\in \mathbb{R}$, then
$V_\rho(H^\varepsilon)(f)(x)=\infty$, a.e. $x\in \mathbb{R}$. In the
next result we establish the behavior of the variation operators
$V_\rho(R_\ell^\varepsilon)$ on $BMO(\mathbb{R}^n)$.

\begin{Th}\label{Th2BMO} Let $\ell=1,\ldots,n$ and $\rho>2$. If $f \in BMO(\mathbb R^n)$ and $V_\rho(R_\ell^\varepsilon)(f)(x) < \infty$, a.e. $x \in \mathbb
 R^n$,  then $V_\rho(R_\ell^\varepsilon)(f) \in BMO(\mathbb R^n)$ and
$\|V_\rho(R_\ell^\varepsilon)f\|_{BMO(\mathbb{R}^n)}\le
C\|f\|_{BMO(\mathbb{R}^n)}.$
\end{Th}

We denote by $\{W_t^\mathcal{L}\}_{t>0}$ the heat semigroup
associated with $\mathcal{L}$. For every $t>0$, we can write
$$
W_t^\mathcal{L}(f)(x)=\int_{\mathbb{R}^n}W_t^\mathcal{L}(x,y)f(y)dy,\,\,\,f\in
L^2(\mathbb{R}^n).
$$
The main properties of the kernel function $W_t^\mathcal{L}(x,y)$,
$x,y\in \mathbb{R}^n$, can be encountered, for instance, in
\cite{DGMTZ}.

The $L^p$-boundedness properties of the variation operator
$V_\rho(W_t^\mathcal{L})$ were studied in \cite{BFHR}.

\begin{Th} \label{SemSch} (\cite[Theorem 1.1]{BFHR}). Let $V\in
RH_q$ where $q>n/2$ and let $\rho>2$. Then, the variation operator
$V_\rho(W_t^\mathcal{L})$ is bounded from $L^p(\mathbb{R}^n)$ into
itself, for every $1<p<\infty$, and from $L^1(\mathbb{R}^n)$ into
$L^{1,\infty}(\mathbb{R}^n)$.
\end{Th}

Our next result shows the behavior of the variation operator $V_\rho
(W_t^\mathcal{L})$ on ${\rm BMO}^\mathcal{L}( \mathbb{R}^n)$.

\begin{Th}\label{Th3BMO} Let $V\in
RH_q$ where $q>n/2$ and let $\rho>2$. Then, the variation operator
$V_\rho (W_t^\mathcal{L})$ is bounded from ${\rm BMO}^\mathcal{L}(
\mathbb{R}^n)$ into itself.
\end{Th}

Let $\ell=1,\ldots,n$. The Riesz transform $R_\ell^\mathcal{L}$ is
defined by
$$
R_\ell^\mathcal{L}(f)=\frac{\partial}{\partial
x_\ell}\mathcal{L}^{-1/2}f,\,\,\,f\in C_c^\infty(\mathbb{R}^n),
$$
where $C_c^\infty(\mathbb{R}^n)$ denotes the space of smooth
functions with compact support in $\mathbb{R}^n$. Here, the negative
square root $\mathcal{L}^{-1/2}$ of $\mathcal{L}$ is defined in
terms of the heat semigroup by
$$
\mathcal{L}^{-1/2}(f)(x)=\frac{1}{\sqrt{\pi}}\int_0^\infty
W_t^\mathcal{L}(f)(x)t^{-1/2}dt.
$$
Fractional powers of the Schr\"odinger operator $\mathcal{L}$ have
been studied in \cite{BHS1}.

In \cite[Proposition 1.1]{BFHR} it was proved that
$R_\ell^\mathcal{L}$ can be extended to $L^p(\mathbb{R})^n$ as the
principal value operator
\begin{equation}\label{limit}
R_\ell^\mathcal{L}(f)(x)=\lim_{\varepsilon\to
0^+}\int_{|x-y|>\varepsilon}R^\mathcal{L}_\ell(x,y)f(y)dy,\,\,\,a.e.\,\,\,x\in
\mathbb{R}^n,
\end{equation}
where
$$
R^\mathcal{L}_\ell(x,y)=-\frac{1}{2\pi}\frac{\partial}{\partial
x_\ell}\int_{\mathbb{R}}(-i\tau)^{-1/2}\Gamma(x,y,\tau)d\tau,\,\,\,x,y\in
\mathbb{R}^n,\,\,x\neq y,
$$
provided that

(i) $1\le p<\infty$, and $V\in RH_n$;

(ii) $1<p<p_0$, where $\frac{1}{p_0}=\frac{1}{q}-\frac{1}{n}$, and
$V\in RH_q$, $n/2\le q<n$.

Moreover, $R_\ell^\mathcal{L}$ is bounded from $L^p(\mathbb{R}^n)$
into itself when $1<p<\infty$ and from $L^1(\mathbb{R}^n)$ into
$L^{1,\infty}(\mathbb{R}^n)$, provided that $V\in RH_n$. Also,
$R_\ell^\mathcal{L}$ is bounded from $L^p(\mathbb{R}^n)$ into itself
when $1<p<p_0$ and $V\in RH_q$, with $n/2\le q<n$ (\cite[Theorems 0.5
and 0.8]{Sh}).

The formal adjoint operator $\mathcal{R}_\ell^\mathcal{L}$ of
$R_\ell^\mathcal{L}$ defined by
$$
\mathcal{R}_\ell^\mathcal{L}(f)(x)=\lim_{\varepsilon\to
0^+}\int_{|x-y|>\varepsilon}R^\mathcal{L}_\ell(y,x)f(y)dy,\,\,\,a.e.\,\,\,x\in
\mathbb{R}^n,
$$
is bounded from $L^p(\mathbb{R}^n)$ into itself when $p_0'<p<\infty$
and $V\in RH_q$, with $n/2\le q<n$, where as usual $p_0'$ denotes
the exponent conjugated to $p_0$.

By defining the truncated Riesz transforms
$R_\ell^{\mathcal{L},\varepsilon}$, $\varepsilon>0$, in the obvious
way, the $L^p$-boundedness properties for the variation operator
$V_\rho(R_\ell^{\mathcal{L},\varepsilon})$ were established in
\cite{BFHR}.

\begin{Th} \label{VarRiesz} (\cite[Theorem 1.2]{BFHR}). Let $\ell=1,\ldots,n$. Assume that $\rho>2$ and $V\in RH_q$, $q\ge n/2$. Then, the variation operator $V_\rho(R^{\mathcal{L},\varepsilon}_\ell)$
is bounded from $L^p(\mathbb{R}^n)$ into itself, when $1<p<p_0$, where
$\frac{1}{p_0}=\Big(\frac{1}{q}-\frac{1}{n}\Big)_+$. Moreover, the operator $V_\rho(R^{\mathcal{L},\varepsilon}_\ell)$
is bounded from $L^1(\mathbb{R}^n)$ into $L^{1,\infty}(\mathbb{R}^n)$.
\end{Th}

If $f\in BMO^\mathcal{L}(\mathbb{R}^n)$ and $\ell=1,\cdots,n$, the
limit in (\ref{limit}) exists for a.e. $x\in \mathbb{R}^n$ (see
\cite{BFHR}). Thus, the Riesz transform $R_\ell^\mathcal{L}$ is
defined by (\ref{limit}) in $BMO^\mathcal{L}(\mathbb{R}^n)$. As it
was mentioned earlier, in the classical case the situation is
different. The behavior of the variation operators associated with
the Riesz transforms $R_\ell^\mathcal{L}$ is established in the
following.

\begin{Th}\label{Th4BMO} Let $\rho >2$, $\ell=1,\ldots,n$ and $V\in RH_q$.

 (i) If $q\geq n$, then the variation operator $V_\rho
(R_\ell^{\mathcal{L},\varepsilon})$ is bounded from
$BMO^\mathcal{L}( \mathbb{R}^n)$ into itself.

(ii) If $q\ge n/2$, then the variation
operator $V_\rho (\mathcal{R}_\ell^{\mathcal{L},\varepsilon})$ is
bounded from $ BMO^\mathcal{L}(\mathbb{R}^n)$ into itself.
\end{Th}

Note that there exists a difference between the results in the
classical and in the Schr\"odinger settings. In the latter case the
operators are defined in the whole $BMO^\mathcal{L}(\mathbb{R}^n)$,
while in the classical case it is necessary to impose a "finite
hypothesis". This fact was observed at first time in
\cite{DGMTZ}.

In order to analyze operators in the Schr\"odinger context on ${\rm
BMO}^\mathcal{L}(\mathbb{R}^n)$ we shall use some ideas developed in
\cite{DGMTZ} and we will again exploit that the Schr\"odinger
operator $\mathcal{L}=-\Delta+V$, where $V\in RH_q$, with $q>n/2$,
is actually a nice perturbation of the Laplacian operator $-\Delta$.


Throughout the proof of the results that we have just showed
the following properties will play an important role.

According to \cite[Proposition 2.14]{DGMTZ} we choose a sequence
$\{x_k\}_{k=1}^\infty \subset \mathbb{R}^n$, such that if
$Q_k=B(x_k,\gamma (x_k))$, $k\in \mathbb{N}$, the following
properties hold:

(i) $\cup_{k=1}^\infty Q_k=\mathbb{R}^n$;

(ii) For every $m\in \mathbb{N}$ there exist $C,\beta >0$ such that, for every $k\in \mathbb{N}$,
$$
\mbox{card }\{l\in \mathbb{N}:2^mQ_l\cap 2^mQ_k\not=\varnothing\}\leq C2^{m\beta }.
$$

According to \cite[p. 346, after Lemma 9]{DGMTZ}, if $H$ is an
operator and $f\in {\rm BMO}^\mathcal{L}(\mathbb{R}^n)$, then $Hf\in {\rm
BMO}^\mathcal{L}(\mathbb{R}^n)$ provided that there exists a
positive constant $C$ such that, for every $k\in\mathbb{N}$,

($i_k$) $\frac{1}{|Q_k|}\int_{Q_k}|H(f)(x)|dx\leq C||f||_{{\rm BMO}^\mathcal{L}(\mathbb{R}^n)}$, and

($ii_k$) $Hf\in {\rm BMO}(Q_k^*)$ and $||Hf||_{{\rm BMO}(Q_k^*)}\leq C||f||_{{\rm BMO}^\mathcal{L}(\mathbb{R}^n)}$. Here ${\rm BMO}(Q_k^*)$ denotes the usual BMO-space on $Q_k^*$.

Moreover, we have that
$$
\|Hf\|_{BMO^\mathcal{L}(\mathbb{R}^n)}\le
M\|f\|_{BMO^\mathcal{L}(\mathbb{R}^n)},
$$
where the constant $M>0$ only depends on the constant $C$.

Note that if the property ($i_k$) above holds for every $k\in
\mathbb{N}$ then, $|H(f)(x)|<\infty$, for almost every $x\in
\mathbb{R}^n$. This fact is different from the one appeared in the
classical Euclidean case (see Theorem 2.2 and Theorem 2.4
\cite{CBS}).

\section{Proof of Theorems \ref{Th1BMO} and \ref{Th2BMO}}

In this section we show our results about the behavior of the variation operator for the classical heat semigroup and Riesz transforms on $BMO(\mathbb{R}^n)$.

\subsection{Proof of Theorem \ref{Th1BMO}}

Assume that $f \in BMO(\mathbb R^n)$ and
$V_\rho(W_t) (f)(x) < \infty$, a.e. $x \in \mathbb R^n$. Let
$B=B(x_0,r_0)$, with $x_0 \in \mathbb R^n$ and $r_0 >0$. We write
$$
f= (f - f_B) \chi_{B^*}+ (f-f_B)\chi_{(B^*)^c} + f_B = f_1+f_2+f_3.
$$
 According to Theorem \ref{VarLp}, we have
\begin{equation} \label{A1}
\int_{\mathbb R^n}\left|V_\rho(W_t)(f_1)(x)\right|^2 dx \leq
C\int_{B^*}\left|f(x)- f_B\right|^2dx \leq C |B|
\|f\|^2_{BMO(\mathbb R^n)}.
\end{equation}
Then $V_\rho(W_t)(f_1)(x) <\infty$, a.e. $x \in \mathbb R^n$.
Moreover, since $\{W_t\}_{t > 0}$ is Markovian,
$V_\rho(W_t)(f_3)=0$. We choose $x_1 \in B(x_0,r_0)$ such that
$V_\rho(W_t)(f_2)(x_1) < \infty$.

If $E_\rho$  denotes the space introduced in Section 1, we can write
\begin{multline}\label{A2}
\frac{1}{|B|}\int_B \left|V_\rho(W_t)(f)(x)- V_\rho(W_t)(f_2)(x_1)\right|dx \\
\shoveleft{\hspace{30mm}= \frac{1}{|B|}\int_B \left| \left\|W_t(f)(x)\right\|_{E_\rho}- \left\|W_t(f_2)(x_1)\right\|_{E_\rho}\right|dx}\\
\shoveleft{\hspace{30mm}\leq \frac{1}{|B|}\int_B  \left\|W_t(f)(x)- W_t(f_2)(x_1)\right\|_{E_\rho}dx }\\
\shoveleft{\hspace{30mm}= \frac{1}{|B|}\int_B  \left\|W_t(f_1)(x) + W_t(f_2)(x) - W_t(f_2)(x_1)\right\|_{E_\rho}dx }\\
\leq \frac{1}{|B|}\int_B  \left\|W_t(f_1)(x)\right\|_{E_\rho}dx +
\frac{1}{|B|} \int_B \left\|W_t(f_2)(x) -
W_t(f_2)(x_1)\right\|_{E_\rho}dx.
\end{multline}
According to (\ref{A1}) we get
\begin{equation}\label{A3}
\frac{1}{|B|}\int_B  \left\|W_t(f_1)(x)\right\|_{E_\rho}dx \leq
\left( \frac{1}{|B|} \int_B \left|V_\rho(W_t)(f_1)(x)\right|^2dx
\right)^\frac{1}{2} \leq C\|f\|_{BMO(\mathbb R^n)}.
\end{equation}
Also, Minkowski inequality and \cite[p. 91]{CMMTV} lead to
\begin{multline}\label{A4}
\frac{1}{|B|} \int_B \left\|W_t(f_2)(x) - W_t(f_2)(x_1)\right\|_{E_\rho}dx \\
\shoveleft{\hspace{25mm}= \frac{1}{|B|} \int_B \left\|\int_{\mathbb R^n}W_t(x,y)f_2(y)dy - \int_{\mathbb R^n}W_t(x_1,y)f_2(y)dy \right\|_{E_\rho}dx} \\
\shoveleft{\hspace{25mm} \leq  C\frac{1}{|B|} \int_B \int_{\mathbb R^n}\left\|W_t(x,y) - W_t(x_1,y)\right\|_{E_\rho}|f_2(y)|dy dx} \\
\shoveleft{\hspace{25mm} \leq \frac{C}{|B|} \int_B \int_{(B^*)^c} \frac{|x-x_1|}{|x-y|^{n+1}}|f(y)-f_B| dy dx}\\
\shoveleft{\hspace{25mm} \leq C\frac{r_0}{|B|} \int_B \int_{(B^*)^c} \frac{1}{(|y-x_0|- |x_0-x|)^{n+1}}|f(y)-f_B| dy dx}\\
\shoveleft{\hspace{25mm} \leq C\frac{r_0}{|B|} \int_B \sum_{k=1}^{\infty}\int_{2^kr_0 \leq |y-x_0| < 2^{k+1}{r_0}} \frac{1}{(|y-x_0|-|x_0-x|)^{n+1}}|f(y)-f_B| dy dx}\\
             \leq C \sum_{k=1}^\infty \frac{1}{2^k}\frac{1}{(2^kr_0)^n}\int_{|y-x_0| < 2^{k+1}r_0}|f(y)- f_B| dy \leq C\|f\|_{BMO(\mathbb R^n)}.\hspace{12mm}
\end{multline}
In the last inequality we have used the following well known property (see
\cite[VIII, Proposition 3.2]{Tor})
$$
\frac{1}{|2^mB|}\int_{2^mB}|f(y)- f_B|dy \leq Cm\|f\|_{BMO(\mathbb R^n)},\;\; m \in \mathbb N.
$$
>From (\ref{A2}), (\ref{A3}) and (\ref{A4}) we conclude that
$$
\frac{1}{|B|}\int_B \left| V_\rho (W_t)(f)(x)-V_\rho(W_t)(f_2)(x_1)\right|dx \leq  C\|f\|_{BMO(\mathbb R^n)}.
$$
Thus, we prove that $V_\rho(W_t)(f)\in BMO(\mathbb R^n)$.

\begin{Rm}\label{Nota1} After a careful read of the proof of Theorem \ref{Th1BMO} we can deduce the following result that will be useful in the proof of Theorem \ref{Th3BMO}.

\begin{Prop}\label{Pr} Let $\rho>2$ and  $\mathcal{A}$ be a set of real decreasing sequences that
converge to zero. Assume that
$$
V_{\rho,\mathcal{A}}(W_t)(f)(x)=\sup_{\{t_j\}_{j\in \mathbb{N}}\in \mathcal{A}}\Big(\sum_{j=1}^\infty |W_{t_j}(f)(x)-W_{t_{j+1}}(f)(x)|^\rho\Big)^{1/\rho}<\infty, \,\,\,a.e.\,\,\,x\in Q,
$$
where $Q$ is a ball in $\mathbb{R}^n$ and $f\in
BMO(\mathbb{R}^n)$. Suppose that $B$ is a ball contained in $Q$. We define
$f_2=(f-f_B)\chi_{(B^*)^c}$ and we choose  $x_1\in B$ such that
$V_{\rho,\mathcal{A}}(W_t)(f_2)(x_1)<\infty$.  Then,
$$
\frac{1}{|B|}\int_B\|W_t(f)(x)-W_t(f_2)(x_1)\|_{E_{\rho,\mathcal{A}}}dx\le C\|f\|_{BMO(\mathbb{R}^n)},
$$
where, for every function $h:(0,\infty)\mapsto \mathbb{C}$, $\|h\|_{E_{\rho,\mathcal{A}}}$ is defined by
$$
\|h\|_{E_{\rho,\mathcal{A}}}=\sup_{\{t_j\}_{j\in \mathbb{N}}\in \mathcal{A}}\Big(\sum_{j=1}^\infty |h(t_j)-h(t_{j+1})|^\rho\Big)^{1/\rho},
$$
and $C>0$ does not depend on $\mathcal{A},$ $f$, nor $B$.
\end{Prop}
\end{Rm}

\subsection{Proof of Theorem \ref{Th2BMO}}

Assume that $f \in BMO(\mathbb R^n)$ and that
$V_\rho(R_\ell^\varepsilon)(f)(x)< \infty$, a.e. $x \in \mathbb
R^n$. To see that $V_\rho(R_\ell^\varepsilon)(f)\in BMO(\mathbb
R^n)$ we extend to $\mathbb R^n$ the procedure developed in the
proof of \cite[Lemma 1.4]{GiTo}.

Let $B=B(x_0,r_0)$ be a ball in $\mathbb R^n$. We decompose $f$, as
usual, by $f=f_1+f_2+f_3$, where $f_1=(f-f_B)\chi_{B^{**}}$,
$f_2=(f-f_B)\chi_{(B^{**})^c}$ and $f_3=f_B$. According to Theorem
\ref{Th2Lp},  we have
\begin{equation}\label{B1}
\int_{\mathbb R^n} |V_\rho(R_\ell^\varepsilon)(f_1)(x)|^2dx \leq
C\int_{B^{**}}|f(x)-f_B|^2dx \leq C |B|\|f\|^2_{BMO(\mathbb R^n)}.
\end{equation}
Then, $V_\rho(R_\ell^\varepsilon)(f_1)(x) < \infty$, a.e. $x \in
\mathbb R^n$. Moreover, $V_\rho(R_\ell^\varepsilon)(f_3)=0$. Then, since $V_\rho(R_\ell^\varepsilon)(f)(x)<\infty$, a.e. $x\in \mathbb{R}^n$, also $V_\rho(R_\ell^\varepsilon)(f_2)(x)<\infty$, a.e. $x\in \mathbb{R}^n$. We
choose $x_1 \in B$ such that $V_\rho(R_\ell^\varepsilon)(f_2)(x_1) <
\infty$.

If $E_\rho$ denotes the space defined in Section 1, by (\ref{B1}) we can write
\begin{multline}\label{B2}
\frac{1}{|B|}\int_B \left|V_\rho(R_\ell^\varepsilon)(f)(x)-V_\rho(R_\ell^\varepsilon)(f_2)(x_1)\right|dx \\
\shoveleft{\hspace{46mm}\leq\frac{1}{|B|}\int_B \left|\|R_\ell^\varepsilon (f)(x)\|_{E_\rho}- \|R_\ell^\varepsilon(f_2)(x_1)\|_{E_\rho}\right|dx} \\
\shoveleft{\hspace{46mm}\leq\frac{1}{|B|}\int_B \|R_\ell^\varepsilon (f)(x) - R_\ell^\varepsilon(f_2)(x_1)\|_{E_\rho}dx}\\
\shoveleft{\hspace{46mm}\leq\frac{1}{|B|}\int_B \|R_\ell^\varepsilon (f_1)(x)\|_{E_\rho}dx+\frac{1}{|B|} \int_B\|R_\ell^\varepsilon(f_2)(x)- R_\ell^\varepsilon(f_2)(x_1)\|_{E_\rho}dx}\\
\leq C\|f\|_{BMO(\mathbb R^n)} + {\frac{1}{|B|}\int_B
\|R_\ell^\varepsilon (f_2)(x) -
R_\ell^\varepsilon(f_2)(x_1)\|_{E_\rho}dx}.
\end{multline}
We define $R_\ell(z)=c_n\frac{z_\ell}{|z|^{n+1}}$,
$z=(z_1,\cdots,z_n) \in \mathbb R^n\setminus \{0\}$. We have that
\begin{equation}\label{B3}
\|R_\ell^\varepsilon(f_2)(x)-R_\ell^\varepsilon(f_2)(x_1)\|_{E_\rho}
\leq A_1(x) + A_2(x),\;\; x \in B,
\end{equation}
where, for every $x \in B$,
$$
A_1(x)=\left\|\int_{|x-y| >
\varepsilon}(R_\ell(x-y)-R_\ell(x_1-y))f_2(y)dy\right\|_{E_\rho}
$$
and
$$
A_2(x) = \Big\|\int_{\mathbb{R}^n}\left(\chi_{\{\varepsilon<
|x-y|\}}(y) - \chi_{\{\varepsilon< |x_1-y|
\}}(y)\right)R_\ell(x_1-y)f_2(y)dy\Big\|_{E_\rho}.
$$
By using Minkowski inequality and well known properties of the
function $R_\ell(z)$ we get
\begin{eqnarray}\label{B4}
A_1(x) & \leq & \int_{\mathbb R^n}|R_\ell(x-y)-R_\ell(x_1-y)||f(y)-f_B|\chi_{(B^{**})^c}(y) dy \nonumber\\
& \leq & C \sum_{k=2}^\infty \int_{2^kr_0 \leq |x_0 -y| \leq 2^{k+1}r_0}\frac{|x-x_1|}{|x-y|^{n+1}}|f(y) - f_B|dy \nonumber \\
& \leq & C \sum_{k=1}^\infty \frac{1}{2^k}\frac{1}{(2^kr_0)^n}\int_{2^{k+1}B}|f(y)-f_B|dy \nonumber \\
& \leq & C \|f \|_{BMO(\mathbb R^n)}, \;\; x \in B.
\end{eqnarray}

In order to analyze $A_2$ we split, for every $j \in \mathbb N$, the
integral there (under $\|.\|_{E_\rho}$) in four terms as follows. Let
$\{\varepsilon_j\}_{j=1}^\infty$ be a real decreasing sequence that
converges to zero. It follows that
\begin{multline}\label{B5}
\int_{\mathbb R^n}\left|\chi_{\{\varepsilon_{j+1} < |x-y| < \varepsilon_j\}}(y) - \chi_{\{\varepsilon_{j+1} < |x_1-y| < \varepsilon_j\}}(y)\right| |R_\ell(x_1-y)||f_2(y)|dy\\
\shoveleft{\hspace{20mm}\leq C\left(\int_{\mathbb R^n} \chi_{\{\varepsilon_{j+1} < |x-y| < \varepsilon_{j+1} +2r_0\}}(y) \chi_{\{\varepsilon_{j+1} < |x-y| < \varepsilon_j\}}(y)|\frac{1}{|x_1-y|^n}|f_2(y)|dy\right.}\\
\shoveleft{\hspace{23mm}+ \int_{\mathbb R^n} \chi_{\{\varepsilon_{j} \ge |x_1-y| < \varepsilon_{j} +2r_0\}}(y)  \chi_{\{\varepsilon_{j+1} < |x-y| < \varepsilon_j\}}(y)\frac{1}{|x_1-y|^n}|f_2(y)|dy}\\
\shoveleft{\hspace{23mm}+ \int_{\mathbb R^n} \chi_{\{\varepsilon_{j+1} < |x_1-y| < \varepsilon_{j+1} +2r_0\}}(y) \chi_{\{\varepsilon_{j+1} < |x_1-y| < \varepsilon_j\}}(y)\frac{1}{|x_1-y|^n}|f_2(y)|dy}\\
\shoveleft{\hspace{23mm}\left.+ \int_{\mathbb R^n} \chi_{\{\varepsilon_{j} < |x-y| < \varepsilon_{j} +2r_0\}}(y) \chi_{\{\varepsilon_{j+1} < |x_1-y| < \varepsilon_j\}}(y)\frac{1}{|x_1-y|^n}|f_2(y)|dy \right)}\\
\shoveleft{\hspace{23mm}= C(A^j_{2,1}(x) + A^j_{2,2}(x) +
A^j_{2,3}(x) + A^j_{2,4}(x)),\hspace{15mm}x\in B\,\,and\,\,j\in
\mathbb{N}.\hspace{42mm}}
\end{multline}

Observe that if $x\in B$, then $A_{2,m}^j(x)=0$, when $m=1,3$,
$j\in \mathbb{N}$ and $r_0\ge \varepsilon_{j+1}$. Also, if $x\in B$,
then $A_{2,m}^j(x)=0$, when $m=2,4$, $j\in \mathbb{N}$ and
$r_0\ge \varepsilon_{j}$.

Let $j\in \mathbb{N}$. Since $2 |x-y| \geq |x_1-y| \geq
\frac{1}{2}|x-y|$, $y \notin B^{**}$ and $x \in B$, H\"{o}lder
inequality leads to
$$
A_{2,1}^j(x) \leq C \left(\int_{\mathbb R^n}\chi_{\{\varepsilon_{j+1} < |x-y| < \varepsilon_j\}}(y) \frac{1}{|x-y|^{ns}}|f_2(y)|^sdy\right)^{\frac{1}{s}}v_{j+1}^{\frac{1}{s'}},\;\; x\in B,
$$
$$
A^j_{2,2}(x) \leq C \left(\int_{\mathbb
R^n}\chi_{\{\max\{\varepsilon_{j+1},
\frac{1}{2}\varepsilon_j\}<|x-y|< \varepsilon_j\}}(y)
\frac{1}{|x-y|^{ns}}|f_2(y)|^sdy\right)^{\frac{1}{s}}v_j^{\frac{1}{s'}},
\;\; x\in B,
$$
$$
A_{2,3}^j(x) \leq C \left(\int_{\mathbb
R^n}\chi_{\{\varepsilon_{j+1} < |x_1-y| < \varepsilon_j\}}(y)
\frac{1}{|x_1-y|^{ns}}|f_2(y)|^sdy\right)^{\frac{1}{s}}v_{j+1}^{\frac{1}{s'}},\;\;
x\in B,
$$
and
$$
A^j_{2,4}(x) \leq C \left(\int_{\mathbb
R^n}\chi_{\{\max\{\varepsilon_{j+1},
\frac{1}{2}\varepsilon_j\}<|x_1-y|< \varepsilon_j\}}(y)
\frac{1}{|x_1-y|^{ns}}|f_2(y)|^sdy\right)^{\frac{1}{s}}v_j^{\frac{1}{s'}},
\;\; x\in B.
$$
Here $1 < s < \infty$, $s'= \frac{s}{s-1}$, and $v_j=(\varepsilon_j+2r_0)^n - \varepsilon_j^n$. Note that $v_j \leq C\max\{r_0,\varepsilon_j\}^{n-1}r_0$, for a certain $C > 0$.

We define the set $\mathcal{A}=\{j\in \mathbb{N}: r_0<\varepsilon _j\}$. We have that
\begin{eqnarray*}
A_{2,1}^j(x)&\leq&C\frac{v_{j+1}^{1/s'}}{\varepsilon _{j+1}^{(n-1)/s'}}\left(\int_{\mathbb{R}^n}\chi _{\{\varepsilon _{j+1}<|x-y|<\varepsilon _j\}}(y)\frac{|f_2(y)|^s}{|x-y|^{n+s-1}}dy\right)^{1/s}\\
&\leq&Cr_0^{1/s'}\left(\int_{\mathbb{R}^n}\chi _{\{\varepsilon _{j+1}<|x-y|<\varepsilon _j\}}(y)\frac{|f_2(y)|^s}{|x-y|^{n+s-1}}dy\right)^{1/s},
\end{eqnarray*}
for every $x\in B$ and $j+1\in \mathcal{A}$. In a similar way we can see that
$$
A_{2,2}^j(x)\leq Cr_0^{1/s'}\left(\int_{\mathbb{R}^n}\chi _{\{\varepsilon _{j+1}<|x-y|<\varepsilon _j\}}(y)\frac{|f_2(y)|^s}{|x-y|^{n+s-1}}dy\right)^{1/s},\quad x\in B\mbox{ and }j\in \mathcal{A},
$$
$$
A_{2,3}^j(x)\leq Cr_0^{1/s'}\left(\int_{\mathbb{R}^n}\chi _{\{\varepsilon _{j+1}<|x_1-y|<\varepsilon _j\}}(y)\frac{|f_2(y)|^s}{|x_1-y|^{n+s-1}}dy\right)^{1/s},\quad x\in B\mbox{ and }j+1\in \mathcal{A},
$$
and
$$
A_{2,4}^j(x)\leq Cr_0^{1/s'}\left(\int_{\mathbb{R}^n}\chi _{\{\varepsilon _{j+1}<|x_1-y|<\varepsilon _j\}}(y)\frac{|f_2(y)|^s}{|x_1-y|^{n+s-1}}dy\right)^{1/s},\quad x\in B\mbox{ and }j\in \mathcal{A}.
$$
Hence, we get
\begin{eqnarray}\label{B6}
\left(\sum _{j=1}^\infty \Big(A_{2,1}^j(x)+A_{2,2}^j(x)\Big)^\rho \right)^{1/\rho}&\leq&C\left(\sum_{j+1\in \mathcal{A}}|A_{2,1}^j(x)|^\rho +\sum_{j\in \mathcal{A}}|A_{2,2}^j(x)|^\rho\right)^{1/\rho}\nonumber\\
&\leq&C\left(\sum_{j\in \mathbb{N}}\left(\int_{\mathbb{R}^n}\chi _{\{\varepsilon _{j+1}<|x-y|<\varepsilon _j\}}(y)\frac{|f_2(y)|^s}{|x-y|^{n+s-1}}dy\right)^{\rho/s}r_0^{\rho /s'}\right)^ {1/\rho}\nonumber\\
&\leq&
C\left(\int_{\mathbb{R}^n}\frac{|f_2(y)|^s}{|x-y|^{n+s-1}}dy\right)^{1/s}r_0^{1/s'}
\nonumber\\
&\leq& C\left(\left(\sum_{k=1}^\infty \frac{1}{(2^kr_0)^n}\int_{|x_0-y|<2^{k+1}r_0}|f(y)-f_B|^sdy\frac{1}{2^{k(s-1)}}\right)^{1/s}\right.\nonumber\\
&\leq&C||f||_{{\rm BMO}(\mathbb{R}^n)}, \quad x\in B.
\end{eqnarray}

In a similar way we get
\begin{equation}\label{B7}
\left(\sum_{j=1}^\infty |A_{2,3}^j(x)+A_{2,4}^j(x)|^\rho
\right)^{1/\rho }\leq C||f||_{{\rm BMO}(\mathbb{R}^n)},\quad x\in B.
\end{equation}

>From (\ref{B5}), (\ref{B6}) and (\ref{B7}) we infer that
\begin{equation}\label{B8}
A_2(x)\leq C||f||_{{\rm BMO}(\mathbb{R}^n)},\quad x\in B.
\end{equation}

Then (\ref{B2}), (\ref{B3}), (\ref{B4}) and (\ref{B8}) imply that
\begin{equation}\label{B9}
\frac{1}{|B|}\int_B|V_\rho (R_\ell^\varepsilon )(f)(x)-V_\rho
(R_\ell^\varepsilon )(f_2)(x_1)|dx\leq C||f||_{{\rm
BMO}(\mathbb{R}^n)}.
\end{equation}
Thus, we prove that $V_\rho (R_\ell^\varepsilon )(f)\in {\rm
BMO}(\mathbb{R}^n)$.

\vspace{10mm}

By proceeding as in the above proof we can establish the following result that will be useful in the sequel.

\begin{Prop}\label{Rtrun}
Let $a>0$ and $\ell=1,...,n$. We define, for every $\varepsilon
>0$ and $f\in L^1_{\rm loc}(\mathbb{R}^n)$,
$$
R_{\ell,a}^\varepsilon (f)(x)=\int_{\varepsilon
<|x-y|<a}\frac{x_\ell-y_\ell}{|x-y|^{n+1}}f(y)dy.
$$
Then, if $\rho >2$, we have that for every $f\in {\rm BMO}(\mathbb{R}^n)$, $V_\rho (R_{\ell,a}^\varepsilon )(f)\in{\rm
BMO}(\mathbb{R}^n)$, and
$$
\|V_\rho (R_{\ell,a}^\varepsilon )(f)\|_{BMO(\mathbb{R}^n)}\le C\|f\|_{BMO(\mathbb{R}^n)},
$$
where $C>0$ does not depend on $f$ nor on $a$.
\end{Prop}

Note that, for every $a>0$, $f\in {\rm BMO}(\mathbb{R}^n)$ and
$\ell=1,\ldots,n$, $V_\rho(R_{\ell,a}^\varepsilon)(f)(x)<\infty$, a.e. $x\in
\mathbb{R}^n$. Indeed, let $a>0$, $f\in {\rm BMO}(\mathbb{R}^n)$ and
$\ell=1,\ldots,n$. Suppose that $m\in \mathbb{N}$. Since
$V_\rho(R_{\ell,a}^\varepsilon)(f)\le
V_\rho(R_{\ell}^\varepsilon)(f)$ and $f\in
L^2_{\rm{loc}}(\mathbb{R}^n)$, according to Theorem \ref{Th2Lp}, we
have that
\begin{eqnarray*}
\int_{B(0,m)}V_\rho(R_{\ell,a}^\varepsilon)(f)(x)dx&\le& \int_{B(0,m)}V_\rho(R_{\ell,a}^\varepsilon)(f\chi_{B(0,m+a)})(x)dx\\
&\le&\int_{B(0,m)}V_\rho(R_{\ell}^\varepsilon)(f\chi_{B(0,m+a)})(x)dx\\
&\le&\sqrt{m}\Big(\int_{B(0,m)}(V_\rho(R_{\ell}^\varepsilon)(f\chi_{B(0,m+a)})(x))^2dx\Big)^{1/2}\\
&\le&\sqrt{m}\Big(\int_{B(0,m+a)}|f(x)|^2dx\Big)^{1/2}<\infty.
\end{eqnarray*}
Hence, $V_\rho(R_{\ell,a}^\varepsilon)(f)(x)<\infty$, a.e. $x\in
B(0,m)$.

\section{Proof of Theorems \ref{Th3BMO} and \ref{Th4BMO}}

In this section we establish the boundedness in $BMO^\mathcal{L}(\mathbb{R}^n)$ of the variation operators for the heat semigroup and Riesz transforms in the Schr\"odinger setting.

\subsection{Proof of Theorem \ref{Th3BMO}}

Assume that $f\in
{\rm BMO}^\mathcal{L}( \mathbb{R}^n)$. Our objective is to show that $V_\rho(W_t^\mathcal{L})(f)$ satisfies the properties $(i_k)$ and $(ii_k)$, for every $k\in \mathbb{N}$ (see Section 2).

Fix $k\in \mathbb{N}$. We now prove $(i_k)$, that is, there exists
$C>0$, independent of $k$, such that
$$
\frac{1}{|Q_k|}\int_{Q_k}|V_\rho(W_t^\mathcal{L})(f)(x)|dx\le C\|f\|_{BMO^\mathcal{L}(\mathbb{R}^n)}.
$$

 We decompose $W_t^\mathcal{L}(f)$ as follows
$$
W_t^\mathcal{L}(f)(x)=H_{k,t}^\mathcal{L}(f)(x)+L_{k,t}^\mathcal{L}(f)(x), \,\,\,x\in \mathbb{R}^n\,\,\,and\,\,\,t>0,
$$
where
$$
H_{k,t}^\mathcal{L}(f)(x)=W_t^\mathcal{L}(f)(x)\chi_{\{t>\gamma(x_k)^2\}}(t),
\,\,\,x\in \mathbb{R}^n\,\,\,and\,\,\,t>0.
$$

It is clear that
\begin{equation}\label{M1}
V_\rho(W_t^\mathcal{L})(f)\le V_\rho(H_{k,t}^\mathcal{L})(f)+V_\rho(L_{k,t}^\mathcal{L})(f).
\end{equation}

Let $\{t_j\}_{j=1}^\infty$ be a real decreasing sequence that
converges to zero. Suppose that $j_k\in \mathbb{N}$ is such that
$t_{j_{k+1}}\le\gamma(x_k)^2< t_{j_k}$. We can write
\begin{eqnarray*}
&&\Big(\sum_{j=1}^\infty |H_{k,t_j}^\mathcal{L}(f)(x)-H_{k,t_{j+1}}^\mathcal{L}(f)(x)|^\rho\Big)^{1/\rho}\\
&\le&\sum_{j=1}^{j_k-1}|W_{t_j}^\mathcal{L}(f)(x)-W_{t_{j+1}}^\mathcal{L}(f)(x)|+|W_{t_{j_k}}^\mathcal{L}(f)(x)|\\
&\le&\sum_{j=1}^{j_k-1}\Big|\int_{t_{j+1}}^{t_j}\frac{\partial}{\partial t}W_t^\mathcal{L}(f)(x)dt\Big|+|W_{t_{j_k}}^\mathcal{L}(f)(x)|\\
&\le& \int_{\gamma(x_k)^2}^\infty\int_{\mathbb{R}^n}\Big|\frac{\partial}{\partial t}W_t^\mathcal{L}(x,y)\Big||f(y)|dydt+\sup_{t\ge \gamma(x_k)^2}|W_t^\mathcal{L}(f)(x)|\\
&=&\Omega_{1,k}(f)(x)+\Omega_{2,k}(f)(x),\,\,\,x\in \mathbb{R}^n.
\end{eqnarray*}
Hence
\begin{equation}\label{M2}
V_\rho(H_{k,t}^\mathcal{L})(f)\le \Omega_{1,k}(f)+\Omega_{2,k}(f).
\end{equation}

According to \cite[(5.4)]{DGMTZ}, we get
\begin{equation}\label{M3}
\frac{1}{|Q_k|}\int_{Q_k}|\Omega_{2,k}(f)(x)|dx\le C\|f\|_{BMO^\mathcal{L}(\mathbb{R}^n)}.
\end{equation}

By using \cite[(2.11)]{DGMTZ}, since $\gamma(x)\sim\gamma(x_k)$, when $x\in Q_k$, we obtain
\begin{eqnarray*}
&&\Omega_{1,k}(f)(x)\le C\int_{\gamma(x_k)^2}^\infty \int_{\mathbb{R}^n}\frac{|f(y)|}{t^{1+n/2}}e^{-c|x-y|^2/t}\Big(1+\frac{t}{\gamma(x_k)^2}\Big)^{-1}dydt\\
&\le& C\int_{\gamma(x_k)^2}^\infty \frac{1}{t^{1+n/2}}\Big(1+\frac{t}{\gamma(x_k)^2}\Big)^{-1}\Big(\int_{|x-y|<\sqrt{t}}+\sum_{m=0}^\infty \int_{\sqrt{t}2^m\le|x-y|<\sqrt{t}2^{m+1}}\Big)|f(y)|\Big(1+\frac{|x-y|}{\sqrt{t}}\Big)^{-1-n}dydt\\
&\le& C\int_{\gamma(x_k)^2}^\infty \frac{1}{t^{1+n/2}}\Big(1+\frac{t}{\gamma(x_k)^2}\Big)^{-1}\sum_{m=0}^\infty \frac{1}{2^{m(1+n)}}\int_{|x-y|<2^m\sqrt{t}}|f(y)|dydt, \,\,\,x\in Q_k.
\end{eqnarray*}

Moreover, by \cite[Lemma 2]{DGMTZ}, it follows that
\begin{eqnarray*}
&&\sum_{m=0}^\infty \frac{1}{2^{m(1+n)}t^{n/2}}\int_{|x-y|<2^m\sqrt{t}}|f(y)|dy\\
&\le& \sum_{m\in \mathbb{N},\,2^m\sqrt{t}\le \gamma(x)}\frac{1}{2^{m(1+n)}t^{n/2}}\int_{|x-y|<2^m\sqrt{t}}|f(y)|dy+\sum_{m\in \mathbb{N},\,2^m\sqrt{t}> \gamma(x)}\frac{1}{2^{m(1+n)}t^{n/2}}\int_{|x-y|<2^m\sqrt{t}}|f(y)|dy\\
&\le&C\|f\|_{BMO^\mathcal{L}(\mathbb{R}^n)}\Big(\sum_{m\in \mathbb{N},\,2^m\sqrt{t}\le \gamma(x)}\frac{1}{2^{m}}\Big(1+ \log\frac{\gamma(x)}{2^m\sqrt{t}}\Big)+\sum_{m\in \mathbb{N}}\frac{1}{2^m}\Big)\\
&\le& C\|f\|_{BMO^\mathcal{L}(\mathbb{R}^n)},\,\,\,t\ge \gamma(x_k)^2\,\,\,and\,\,\,x\in Q_k.
\end{eqnarray*}

Then,
\begin{eqnarray*}
\Omega_{1,k}(f)(x)&\le& C\int_{\gamma(x_k)^2}^\infty \Big(1+\frac{t}{\gamma(x_k)^2}\Big)^{-1}\frac{dt}{t}\|f\|_{BMO^\mathcal{L}(\mathbb{R}^n)}\\
&\le& C\int_1^\infty \frac{du}{u(1+u)}\|f\|_{BMO^\mathcal{L}(\mathbb{R}^n)}\\
&\le& C \|f\|_{BMO^\mathcal{L}(\mathbb{R}^n)},\,\,\,x\in Q_k.
\end{eqnarray*}

Hence, we obtain
\begin{equation}\label{M4}
\frac{1}{|Q_k|}\int_{Q_k}\Omega_{1,k}(f)(x)dx\le C \|f\|_{BMO^\mathcal{L}(\mathbb{R}^n)}.
\end{equation}

By combining (\ref{M2}), (\ref{M3}) and (\ref{M4}) we get
\begin{equation}\label{M5}
\frac{1}{|Q_k|}\int_{Q_k}V_\rho(H_{k,t}^\mathcal{L})(f)(x)dx\le C\|f\|_{BMO^\mathcal{L}(\mathbb{R}^n)}.
\end{equation}
Here $C>0$ does not depend on $k$.

We now decompose $f$ as follows
$$
f=f\chi_{Q_k^*}+f\chi_{(Q_k^*)^c}=f_1+f_2.
$$

It is clear that
\begin{equation}\label{M6}
V_\rho(L_{k,t}^\mathcal{L})(f)\le V_\rho(L_{k,t}^\mathcal{L})(f_1)+V_\rho(L_{k,t}^\mathcal{L})(f_2).
\end{equation}

By proceeding as above we get
\begin{eqnarray*}
V_\rho(L_{k,t}^\mathcal{L})(f_1)(x)&\le& C\Big(\sup_{\{t_j\}_{j\in \mathbb{N}}\downarrow 0}\Big(\sum_{j\in \mathbb{N},\,t_j\le\gamma(x_k)^2}|W_{t_j}^\mathcal{L}(f_1)(x)-W_{t_{j+1}}^\mathcal{L}(f_1)(x)|^\rho\Big)^{1/\rho}\\
&+& \sup_{0<t\le \gamma(x_k)^2}|W_t^\mathcal{L}(f_1)(x)|\Big)\\
&\le& C(V_\rho(W_t^\mathcal{L})(f_1)(x)+W_*^\mathcal{L}(f_1)(x)),\,\,\,x\in \mathbb{R}^n.
\end{eqnarray*}

Since $W_*^\mathcal{L}$ and $V_\rho(W_t^\mathcal{L})$ are bounded
operators from $L^2(\mathbb{R}^n)$ into itself (see Theorem
\ref{VarLp}) it follows that
\begin{eqnarray*}
\frac{1}{|Q_k|}\int_{Q_k}V_\rho(L_{k,t}^\mathcal{L})(f_1)(x)dx&\le& \Big(\frac{1}{|Q_k|}\int_{\mathbb{R}^n}\Big(V_\rho(L_{k,t}^\mathcal{L})(f_1)(x)\Big)^2dx\Big)^{1/2}\\
&\le& C\Big(\frac{1}{|Q_k|}\int_{Q_k^*}|f(x)|^2dx\Big)^{1/2}.
\end{eqnarray*}
Then, from \cite[Corollary 3]{DGMTZ}, we deduce that
\begin{equation}\label{M7}
\frac{1}{|Q_k|}\int_{Q_k}V_\rho(L_{k,t}^\mathcal{L})(f_1)(x)dx\le C\|f\|_{BMO^\mathcal{L}(\mathbb{R}^n)}.
\end{equation}

On the other hand, we can write
\begin{eqnarray}\label{M8}
V_\rho(L_{k,t}^\mathcal{L})(f_2)(x)&\le&\int_0^{\gamma(x_k)^2}\int_{(Q_k^*)^c}\Big|\frac{\partial}{\partial t}W_t^\mathcal{L}(x,y)\Big||f(y)|dydt+\sup_{0<t\le \gamma(x_k)^2}|W_t^\mathcal{L}(f_2)(x)|\nonumber\\
&=&\Omega_{3,k}(f)(x)+\Omega_{4,k}(f)(x),\,\,\, x\in \mathbb{R}^n.
\end{eqnarray}

According  to \cite[(2.7)]{DGMTZ}, for certain $C,c>0$, we get
\begin{multline}\label{M9}
\int_0^{\gamma(x_k)^2}\int_{(Q_k^*)^c}\Big|\frac{\partial}{\partial t}W_t^\mathcal{L}(x,y)\Big||f(y)|dydt\\
\shoveleft{\hspace{16mm}\leq \int_0^{\gamma(x_k)^2}\int_{|x-y|>\gamma(x_k)}\left|\frac{\partial}{\partial t}W_t^\mathcal{L}(x,y)\right|dt|f(y)|dy}\\
\shoveleft{\hspace{16mm}\leq C\int_0^{\gamma(x_k)^2}\int_{|x-y|>\gamma(x_k)}|f(y)|\frac{e^{-c|x-y|^2/t}}{t^{n/2+1}}dtdy} \\
\shoveleft{\hspace{16mm} \leq C\int_0^{\gamma(x_k)^2}\frac{1}{t^{n/2+1}}\sum_ {j=0}^\infty \int_{2^j\gamma(x_k)<|x-y|\leq 2^{j+1}\gamma(x_k)}|f(y)|\left(\frac{t}{|x-y|^2}\right)^{(n+1)/2}dydt} \\
\shoveleft{ \hspace{16mm}\leq C\int_0^{\gamma(x_k)^2}\frac{1}{\sqrt{t}}\sum_ {j=0}^\infty \frac{1}{(2^j\gamma(x_k))^{n+1}}\int_{|x-y|\leq 2^{j+1}\gamma(x_k)}|f(y)|dydt }\\
 \leq C\sum_ {j=0}^\infty \frac{1}{2^j(2^j\gamma(x_k))^n}\int_{|x_k-y|\leq 2^{j+3}\gamma(x_k)}|f(y)|dy\leq \dot{C||f||_{{\rm BMO}^\mathcal{L}(\mathbb{R}^n)}},\quad x\in
 Q_k.
\end{multline}

Then, by (\ref{M9}),
$$
\frac{1}{|Q_k|}\int_{Q_k}\Omega_{3,k}(f)(x)dx\le C||f||_{{\rm BMO}^\mathcal{L}(\mathbb{R}^n)}.
$$

Moreover, since $|f|\in {\rm BMO}^\mathcal{L}(\mathbb{R}^n)$, \cite[Theorem 6]{DGMTZ} implies that
$$
\frac{1}{|Q_k|}\int_{Q_k}\Omega_{4,k}(f)(x)dx\le \frac{1}{|Q_k|}\int_{Q_k}W_*^\mathcal{L}(|f|)(x)dx\le C||f||_{{\rm BMO}^\mathcal{L}(\mathbb{R}^n)}.
$$

Hence, we conclude that
\begin{equation}\label{M10}
\frac{1}{|Q_k|}\int_{Q_k}V_\rho(L_{k,t}^\mathcal{L})(f_2)(x)dx\le C||f||_{{\rm BMO}^\mathcal{L}(\mathbb{R}^n)}.
\end{equation}

By combining (\ref{M7}) and (\ref{M10}) we deduce
\begin{equation}\label{M11}
\frac{1}{|Q_k|}\int_{Q_k}V_\rho(L_{k,t}^\mathcal{L})(f)(x)dx\le C||f||_{{\rm BMO}^\mathcal{L}(\mathbb{R}^n)}.
\end{equation}

Finally, (\ref{M1}), (\ref{M5}) and (\ref{M11}) imply that
$$
\frac{1}{|Q_k|}\int_{Q_k}V_\rho(W_t^\mathcal{L})(f)(x)dx\le C||f||_{{\rm BMO}^\mathcal{L}(\mathbb{R}^n)}.
$$
Note that $C>0$ does not depend on $k$.

Thus the property $(i_k)$ is established.

We are now going to prove $(ii_k)$. Assume that $B=B(x_0,r_0)\subset
Q_k^*$, with $x_0\in \mathbb{R}^n$ and $r_0>0$. Our objective is to
see that
\begin{equation}\label{C1}
\frac{1}{|B|}\int_B|V_\rho(W_t^\mathcal{L})(f)(x)-c_B|dx\le C||f||_{{\rm BMO}^\mathcal{L}(\mathbb{R}^n)},
\end{equation}
for a certain constant $c_B$, where $C>0$ does not depend on $k$ nor on
$B$. We decompose $W_t^\mathcal{L}$ as follows
\begin{equation}\label{C1.1}
W_t^\mathcal{L}(f)=H_{k,t}^\mathcal{L}(f)+(L_{k,t}^\mathcal{L}(f)-L_{k,t}(f))+L_{k,t}(f),\,\,\,t>0,
\end{equation}
where $H_{k,t}^\mathcal{L}$ and $L_{k,t}^\mathcal{L}$ are defined as before, and
$$
L_{k,t}(f)=W_t(f)\chi_{\{0<t\le \gamma(x_k)^2\}}(t), \,\,\,t>0.
$$
Suppose that $c_B=\|h_B\|_{E_\rho}$, where $h_B:(0,\infty)\mapsto \mathbb{C}$ is a function that will be specified later. Then, we can write
\begin{eqnarray*}
|V_\rho(W_t^\mathcal{L})(f)(x)-c_B|&=&|\|W_t^\mathcal{L}(f)(x)\|_{E_\rho}-\|h_B\|_{E_\rho}|\\
&\le&\|W_t^\mathcal{L}(f)(x)-h_B(t)\|_{E_\rho}\\
&\le&\|H_{k,t}^\mathcal{L}(f)(x)+(L_{k,t}^\mathcal{L}(f)(x)-L_{k,t}(f)(x))+L_{k,t}(f)(x)-h_B(t)\|_{E_\rho}\\
&\le&\|H_{k,t}^\mathcal{L}(f)(x)\|_{E_\rho}+\|L_{k,t}^\mathcal{L}(f)(x)-L_{k,t}(f)(x)\|_{E_\rho}+\|L_{k,t}(f)(x)-h_B(t)\|_{E_\rho}.
\end{eqnarray*}
Then, (\ref{C1}) will be shown when we prove the following three properties:

(A1) $\frac{1}{|B|}\int_B\|H_{k,t}^\mathcal{L}(f)(x)\|_{E_\rho}dx\le C||f||_{{\rm BMO}^\mathcal{L}(\mathbb{R}^n)};$

(A2) $\frac{1}{|B|}\int_B\|L_{k,t}^\mathcal{L}(f)(x)-L_{k,t}(f)(x)\|_{E_\rho}dx\le C||f||_{{\rm BMO}^\mathcal{L}(\mathbb{R}^n)};$ and

(A3) $\frac{1}{|B|}\int_B\|L_{k,t}(f)(x)-h_B(t)\|_{E_\rho}dx\le C||f||_{{\rm BMO}^\mathcal{L}(\mathbb{R}^n)},$

\noindent for a certain function $h_B:(0,\infty)\mapsto \mathbb{C}$, and where $C>0$ is independent of $k$ and $B$.

According to (\ref{M2}) we have
$$
V_\rho(H_{k,t}^\mathcal{L})(f)\le \Omega_{1,k}(f)+\Omega_{2,k}(f).
$$
By proceeding as above we get
\begin{equation}\label{C2}
|\Omega_{1,k}(f)(x)|\le C||f||_{{\rm BMO}^\mathcal{L}(\mathbb{R}^n)},\,\,\,x\in Q_k^*.
\end{equation}
Moreover, by \cite[(5.4)]{DGMTZ},
\begin{equation}\label{C3}
|\Omega_{2,k}(f)(x)|\le C||f||_{{\rm BMO}^\mathcal{L}(\mathbb{R}^n)},\,\,\,x\in Q_k^*.
\end{equation}
>From (\ref{C2}) and (\ref{C3}) it follows that $(A1)$ holds.

To establish $(A2)$ we firstly observe that
\begin{eqnarray*}
V_\rho(L_{k,t}^\mathcal{L}-L_{k,t})(f)(x)&\le&\int_0^{\gamma(x_k)^2}\int_{\mathbb{R}^n}\Big|\frac{\partial}{\partial t}(W_t^\mathcal{L}(x,y)-W_t(x,y))\Big||f(y)|dydt\\
&+&\sup_{0<t\le \gamma(x_k)^2}|W_t^\mathcal{L}(f)(x)-W_t(f)(x)|\\
&=& \Omega_{5,k}(f)(x)+\Omega_{6,k}(f)(x),\,\,\,x\in \mathbb{R}^n.
\end{eqnarray*}

By \cite[(5.5)]{DGMTZ} we get
\begin{equation}\label{C4}
\Omega_{6,k}(f)(x)\le C||f||_{{\rm BMO}^\mathcal{L}(\mathbb{R}^n)},\,\,\,x\in Q_k^*.
\end{equation}

The perturbation formula
(\cite[(5.25)]{DGMTZ}) allows us to write
\begin{eqnarray*}
\frac{\partial}{\partial t}\Big(W_t(x,y)-W_t^{\mathcal{L}}(x,y)\Big)&=&\int_{\mathbb{R}^n}V(z)W_{t/2}^{\mathcal{L}}(x,z)W_{t/2}(z,y)dz\\
&+&\int_0^{t/2}\int_{\mathbb{R}^n}V(z)\frac{\partial}{\partial t}W_{t-s}^{\mathcal{L}}(x,z)W_{s}(z,y)dzds\\
&+&\int_{t/2}^t\int_{\mathbb{R}^n}V(z)W_{t-s}^{\mathcal{L}}(x,z)\frac{\partial}{\partial s}W_{s}(z,y)dzds\\
&=&\sum_{j=1}^3K_j(x,y,t),\,\,\,x,y\in \mathbb{R}^n\,\,\,and\,\,\,t>0.
\end{eqnarray*}
According to \cite[(2.3) and (2.8)]{DGMTZ}, we get
\begin{eqnarray*}
|K_1(x,y,t)|&\le& Ct^{-n}\int_{\mathbb{R}^n}V(z)e^{-\frac{|x-z|^2+|z-y|^2}{4t}}dz\\
&\le&Ct^{-n/2}e^{-\frac{|x-y|^2}{8t}} \int_{\mathbb{R}^n}V(z)t^{-n/2}e^{-\frac{|x-z|^2}{8t}}dz\\
&\le&C\gamma(x)^{-\delta}t^{-1+(\delta-n)/2}e^{-\frac{|x-y|^2}{8t}},\,\,\,x,y\in
\mathbb{R}^n\,\,\,and\,\,\,0<t<\gamma(x_k)^2.
\end{eqnarray*}
Here and in the sequel $\delta$ represents a positive constant.

Moreover, by using \cite[(2.7) and (2.8)]{DGMTZ}, since $t/2<t-s<t$ when $0<s<t/2$, it follows that
\begin{eqnarray*}
|K_2(x,y,t)|&\le& C\int_0^{t/2}\int_{\mathbb{R}^n}V(z)\frac{1}{(t-s)^{1+n/2}}e^{-c\frac{|x-z|^2}{t-s}}\frac{1}{s^{n/2}}e^{-\frac{|z-y|^2}{4s}}dzds\\
&\le&C\int_0^{t/2}\int_{\mathbb{R}^n}V(z)\frac{1}{t^{1+n/2}}e^{-c\frac{|x-z|^2}{t}}\frac{1}{s^{n/2}}e^{-\frac{|z-y|^2}{4s}}dzds\\
&\le&C\frac{1}{t^{1+n/2}}e^{-c\frac{|x-y|^2}{t}}\int_0^{t/2}\int_{\mathbb{R}^n}V(z)\frac{1}{s^{n/2}}e^{-c\frac{|z-y|^2}{s}}dzds\\
&\le&C\frac{1}{t^{1+n/2}}e^{-c\frac{|x-y|^2}{t}}\int_0^{t/2}\frac{s^{-1+\delta/2}}{\gamma(y)^\delta}ds\\
&\le&\gamma(y)^{-\delta}t^{-1+(\delta-n)/2}e^{-c\frac{|x-y|^2}{t}},\,\,\,x,y\in
\mathbb{R}^n\,\,\,and\,\,\,0<t<\gamma(x_k)^2.
\end{eqnarray*}

By proceeding in a similar way we obtain
$$
|K_3(x,y,t)|\le
C\gamma(y)^{-\delta}t^{-1+(\delta-n)/2}e^{-c\frac{|x-y|^2}{t}},\,\,\,x,y\in
\mathbb{R}^n\,\,\,and\,\,\,0<t<\gamma(x_k)^2.
$$
Hence, since $\gamma(x)\sim\gamma(y)\sim\gamma(x_k)$, provided that
$|x-y|\le\gamma(x_k)$ and $x\in Q_k^*$, we conclude that
$$
\Big|\frac{\partial}{\partial
t}\Big(W_t^{\mathcal{L}}(x,y)-W_t(x,y)\Big)\Big|\le
C\gamma(y)^{-\delta}t^{-1+(\delta-n)/2}e^{-c\frac{|x-y|^2}{t}},\,\,\,x\in
Q_k^*,\,\,|x-y|\le \gamma(x_k)\,\,\,and\,\,\,0<t<\gamma(x_k)^2.
$$
We obtain that, for certain $C,c>0$
\begin{eqnarray*}
\Omega_{5,k}(f)(x)&\leq&C\int_0^{\gamma(x_k)^2}\frac{t^{\delta /2-1}}{\gamma(x_k)^\delta}\int_{\mathbb{R}^n}\frac{e^{-c|x-y|^2/t}}{t^{n/2}}|f(y)|dydt\\
&\leq &C\int_0^{\gamma(x_k)^2}\frac{t^{\delta /2-1}}{\gamma(x_k)^\delta}
\sum_{j=0}^\infty\frac{e^{-c2^{2j}}}{t^{n/2}}\int_{|x-y|\leq 2^j\sqrt{t}}|f(y)|dydt\\
&\leq &C\int_0^{\gamma(x_k)^2}\frac{t^{\delta /2-1}}{\gamma(x_k)^\delta}
\sum_{j=0}^\infty \frac{2^{jn}e^{-c2^{2j}}}{(2^j\sqrt{t})^n}\int_{|x-y|\leq 2^j\sqrt{t}}|f(y)|dydt,\quad x\in Q_k^*.
\end{eqnarray*}
Moreover, by \cite[Lemma 3.14]{DGMTZ}, since $\gamma(x)\sim\gamma(x_k)$, $x\in Q_k^*$,
\begin{multline*}
\sum_{j=0}^\infty \frac{2^{jn}e^{-c2^{2j}}}{(2^j\sqrt{t})^n}\int_{|x-y|\leq 2^j\sqrt{t}}|f(y)|dy\\
\shoveleft{\hspace{10mm}=\sum_{j\in \mathbb{N},\,2^j\sqrt{t}\leq \gamma(x)}\frac{2^{jn}e^{-c2^{2j}}}{(2^j\sqrt{t})^n}\int_{|x-y|\leq 2^j\sqrt{t}}|f(y)|dy
+\sum_{j\in \mathbb{N},\,2^j\sqrt{t}>\gamma(x)}\frac{2^{jn}e^{-c2^{2j}}}{(2^j\sqrt{t})^n}\int_{|x-y|\leq 2^j\sqrt{t}}|f(y)|dy}\\
\shoveleft{\hspace{10mm}\leq C||f||_{{\rm BMO}^\mathcal{L}(\mathbb{R}^n)}\left(\sum_{j\in \mathbb{N},\,2^j\sqrt{t}\leq \gamma(x)}2^{jn}e^{-c2^{2j}}\left(1+\log \frac{\gamma(x)}{2^j\sqrt{t}}\right)+
\sum_{j\in \mathbb{N},\,2^j\sqrt{t}>\gamma(x)}2^{jn}e^{-c2^{2j}}\right)}\\
\shoveleft{\hspace{10mm}\leq C||f||_{{\rm
BMO}^\mathcal{L}(\mathbb{R}^n)}\left(\frac{\gamma(x_k)}{\sqrt{t}}\right)^\varepsilon
,\,\,\,x\in Q_k^*\,\,\,and\,\,\,0<t<\gamma(x_k)^2,}\hspace{20mm}
\end{multline*}
where $\varepsilon \in (0,\delta )$.

Therefore, we have that
\begin{equation}\label{C5}
\Omega_{5,k}(f)(x)\leq C||f||_{{\rm BMO}^\mathcal{L}(\mathbb{R}^n)}\int_0^{\gamma(x_k)^2}\frac{t^{\delta /2-1-\varepsilon /2}}{\gamma(x_k)^{\delta -\varepsilon}}dt\leq C||f||_{{\rm BMO}^\mathcal{L}(\mathbb{R}^n)},\quad x\in Q_k^*.
\end{equation}

>From (\ref{C4}) and (\ref{C5}) we infer $(A2)$.

By (\ref{C1.1}) it follows that
$$
V_\rho(L_{k,t})(f)\le V_\rho(W_t^\mathcal{L})(f)+V_\rho(H_{k,t}^\mathcal{L})(f)+V_\rho(L_{k,t}^\mathcal{L}-L_{k,t})(f).
$$
By proceeding as in the proof of $(i_k)$ we get
$$
\int_{Q_k^*}V_\rho(W_t^\mathcal{L})(f)(x)dx<\infty.
$$
Then, $V_\rho(W_t^\mathcal{L})(f)(x)<\infty$, a.e. $x\in Q_k^*$.

>From (\ref{C4}) and (\ref{C5}) we deduce $V_\rho(L_{k,t}^\mathcal{L}-L_{k,t})(f)(x)<\infty$, a.e. $x\in Q_k^*$. Also, by (\ref{C2}) and (\ref{C3}), $V_\rho(H_{k,t}^\mathcal{L})(f)(x)<\infty$, a.e. $x\in Q_k^*$.

Hence, $V_\rho(L_{k,t})(f)(x)<\infty$, a.e. $x\in Q_k^*$. We consider the following decomposition of $f$
$$
f=(f-f_B)\chi_{B^*}+(f-f_B)\chi_{(B^*)^c}+f_B=f_1+f_2+f_3.
$$

Note that
\begin{eqnarray*}
V_\rho(L_{k,t})(f_1)(x)&\le& C\Big(\sup_{\{t_j\}_{j\in \mathbb{N}},\,t_j\le\gamma(x_k)^2}\Big(\sum_{j=1}^\infty |W_{t_j}(f_1)(x)-W_{t_{j+1}}(f_1)(x)|^\rho\Big)^{1/\rho}\\
&+&\sup_{0<t\le \gamma(x_k)^2}|W_t(f_1)(x)|\Big)\\
&\le& C(V_\rho(W_t)(f_1)(x)+W_*(f_1)(x)).
\end{eqnarray*}

Then, since $W_*$ and $V_\rho(W_t)$ are bounded operators from $L^2(\mathbb{R}^n)$ into itself (see Theorem \ref{VarLp}), we obtain
\begin{eqnarray*}
\int_{Q_k^*}|V_\rho(L_{k,t})(f_1)(x)|dx&\le&C\Big(|Q_k|\int_{B^*}|f(x)-f_B|^2dx\Big)^{1/2}\\
&\le&C(|B||Q_k|)^{1/2}\|f\|_{BMO^\mathcal{L}(\mathbb{R}^n)}<\infty.
\end{eqnarray*}
Hence $V_\rho(L_{k,t})(f_1)(x)<\infty$, a.e. $x\in Q_k^*$.

Also, since $\int_{\mathbb{R}^n}W_t(x,y)dy=1$, $x\in \mathbb{R}^n$ and $t>0$, we get
$$
V_\rho(L_{k,t})(f_3)=|f_B|<\infty, \,\,\,x\in \mathbb{R}^n.
$$

Therefore, we deduce that $V_\rho(L_{k,t})(f_2)(x)<\infty$, a.e. $x\in Q_k^*$.

We choose $z_1\in B$ such that $V_\rho(L_{k,t})(f_2)(z_1)< \infty$ and we define $h_B(t)=L_{k,t}(f_2)(z_1)$, $t\in (0,\infty)$.

Suppose that $\{t_j\}_{j\in \mathbb{N}}$ is a real decreasing
sequence that converges to zero and that $j_k\in \mathbb{N}$ is such
that $t_{j_k}\le \gamma(x_k)^2$ and $t_{j_k-1}>\gamma(x_k)^2$. We
can write
\begin{eqnarray*}
&&\Big(\sum_{j=1}^\infty
\Big|L_{k,t_j}(f)(x)-L_{k,t_j}(f_2)(z_1)-(L_{k,t_{j+1}}(f)(x)-L_{k,t_{j+1}}(f_2)(z_1))\Big|^\rho\Big)^{1/\rho}\\
&=&\Big(\sum_{j=j_k}^\infty
\Big|W_{t_j}(f)(x)-W_{t_j}(f_2)(z_1)-(W_{t_{j+1}}(f)(x)-W_{t_{j+1}}(f_2)(z_1))\Big|^\rho+
\Big|W_{t_{j_k}}(f)(x)-W_{t_{j_k}}(f_2)(z_1)\Big|^\rho\Big)^{1/\rho}\\
&\le& C\Big(\Big(\sum_{j=j_k}^\infty
\Big|W_{t_j}(f)(x)-W_{t_j}(f_2)(z_1)-(W_{t_{j+1}}(f)(x)-W_{t_{j+1}}(f_2)(z_1))\Big|^\rho\Big)^{1/\rho}\\
&+&\sup_{0<t\le \rho(x_k)^2}|W_t(f)(x)-W_t(f_2)(z_1)|.
\end{eqnarray*}

We have that
\begin{eqnarray*}
&&\Big|\|L_{k,t}(f)(x)\|_{E_\rho}-\|h_B\|_{E_\rho}\Big|\\
&\le&\sup_{\{t_j\}_{j\in
\mathbb{N}}\downarrow 0,\,0<t_j\le \gamma(x_k)^2}\Big(\sum_{j=1}^\infty
\Big|W_{t_j}(f)(x)-W_{t_{j}}(f_2)(z_1)-(W_{t_{j+1}}(f)(x)-W_{t_{j+1}}(f_2)(z_1))\Big|^\rho\Big)^{1/\rho}\\
&+&\sup_{0<t\le \rho(x_k)^2}|W_t(f)(x)-W_t(f_2)(z_1)|.
\end{eqnarray*}

By taking into account Proposition \ref{Pr} we obtain
\begin{eqnarray}\label{C6}
\frac{1}{|B|}\int_B\sup_{\{t_j\}_{j\in \mathbb{N}}\downarrow
0,\,0<t_j\le \gamma(x_k)^2}\Big(\sum_{j=1}^\infty
&&\Big|W_{t_j}(f)(x)-W_{t_{j}}(f_2)(z_1)-(W_{t_{j+1}}(f)(x)-W_{t_{j+1}}(f_2)(z_1))\Big|^\rho\Big)^{1/\rho}dx\nonumber\\
&\le&
C\|f\|_{BMO^\mathcal{L}(\mathbb{R}^n)}.
\end{eqnarray}

 Also, according to \cite[Pages 348 and 349]{DGMTZ} it follows that
\begin{equation}\label{C7}
\frac{1}{|B|}\int_B\sup_{0<t\le
\rho(x_k)^2}|W_t(f)(x)-W_t(f_2)(z_1)|dx\le
C\|f\|_{BMO^\mathcal{L}(\mathbb{R}^n)}.
\end{equation}

>From (\ref{C6}) and (\ref{C7}) we deduce $(A3)$.

Note that in all the occurrences the constant $C>0$ does not depend on $k$ nor on $B$.

Thus the proof of $(ii_k)$ is finished.

\subsection{Proof of Theorem \ref{Th4BMO}}

 $(i)$ Assume, without lost of generality (\cite{G}), that $V\in RH_q$ where $q>n$.  We are going to prove that the variation operator $V_\rho(R_\ell^{\mathcal{L},\varepsilon})$ is bounded from $BMO^\mathcal{L}(\mathbb{R}^n)$ into itself. We consider, for every $k\in \mathbb{N}$, the local
operators
$$
R_{\ell,k}^\mathcal{L}(f)(x)=PV\int_{|x-y|<\gamma(x_k)}R_\ell^\mathcal{L}(x,y)f(y)dy,
$$
and
$$
R_{\ell,k}(f)(x)=PV\int_{|x-y|<\gamma(x_k)}R_\ell(x-y)f(y)dy.
$$
Note that $|y-x_k|\leq 3\gamma(x_k)$ when $x\in Q_k^*$ and
$|x-y|<\gamma(x_k)$. Then, if $f\in{\rm
BMO}^\mathcal{L}(\mathbb{R}^n)$,
$$
R_{\ell,k}(f)(x)=\lim_{\varepsilon \rightarrow 0^+}\int_{\varepsilon
<|x-y|<\gamma(x_k)} R_\ell(x-y)f(y)\chi
_{3Q_k}(y)dy,\,\,\,a.e.\,\,\,x\in Q_k^*,
$$
that is,  this limit exists for almost all $x\in Q_k^*$ when $f\in
{\rm BMO}^\mathcal{L}(\mathbb{R}^n)$. Also,
$R_{\ell,k}^\mathcal{L}(f)(x)$ is defined for almost every $x\in
Q_k^*$ when $f\in {\rm BMO}^\mathcal{L}(\mathbb{R}^n)$ (see
\cite[Proposition 1.1]{BFHR}).

Let $f\in {\rm BMO}^\mathcal{L}(\mathbb{R}^n)$. We are going to
show, for every $k\in \mathbb{N}$, the properties $(i_k)$ and $(ii_k)$ when $H=V_\rho
(R_\ell^{\mathcal{L},\varepsilon })$. Let $k\in \mathbb{N}$. We can
write
\begin{eqnarray*}
V_\rho (R_\ell^{\mathcal{L},\varepsilon})(f)&=&(V_\rho (R_\ell^{\mathcal{L},\varepsilon})(f)-V_\rho (R_{\ell,k}^{\mathcal{L},\varepsilon})(f))\\
&+&(V_\rho (R_{\ell,k}^{\mathcal{L},\varepsilon})(f)-V_\rho (R_{\ell,k}^\varepsilon) (f))+V_\rho (R_{\ell,k}^\varepsilon)(f)\\
&=&F_{1,k}+F_{2,k}+V_\rho (R_{\ell,k}^\varepsilon )(f).
\end{eqnarray*}
We have that
\begin{eqnarray*}
|F_{1,k}(x)|&\leq&V_\rho (R_\ell^{\mathcal{L},\varepsilon}-R_{\ell,k}^{\mathcal{L},\varepsilon})(f)(x)\\
&=&\sup_{\{\varepsilon _j\}_{j\in \mathbb{N}}\downarrow 0}\left(\sum_{j=1}^\infty \left|\int_{\varepsilon _{j+1}<|x-y|<\varepsilon _j}R_\ell^\mathcal{L}(x,y)f(y)dy-\int_{\varepsilon _{j+1}<|x-y|<\varepsilon _j,\,|x-y|<\gamma(x_k)}R_\ell^\mathcal{L}(x,y)f(y)dy\right|^\rho \right)^{1/\rho }\\
&=&\sup_{\{\varepsilon _j\}_{j\in \mathbb{N}}\downarrow 0}\left(\sum_{j=1}^\infty \left|\int_{\varepsilon _{j+1}<|x-y|<\varepsilon _j,|x-y|\geq \gamma(x_k)}R_\ell^\mathcal{L}(x,y)f(y)dy\right|^\rho \right)^{1/\rho }\\
&\leq&\int_{|x-y|>\gamma(x_k)}|R_\ell^\mathcal{L}(x,y)||f(y)|dy,\quad
x\in Q_k^*.
\end{eqnarray*}
Then, according to \cite[Lemma 3, (a)]{BHS2}, since $\gamma(x_k)\geq
M\gamma(x)$, $x\in Q_k^*$, for a certain $0<M<1$ that does not depend
on $k\in\mathbb{N}$, it follows that
\begin{eqnarray*}
|F_{1,k}(x)|&\leq&C\int_{|x-y|>M\gamma(x)}\frac{1}{|x-y|^n}\frac{1}{1+|x-y|/\gamma(x)}|f(y)|dy\\
&\leq&C\sum_{j=0}^\infty \int_{M2^j\gamma(x)<|x-y|<M2^{j+1}\gamma(x)}
\frac{1}{|x-y|^n}\frac{1}{1+|x-y|/\gamma(x)}|f(y)|dy\\
&\leq&C\sum_{j=0}^\infty \frac{1}{2^j}\frac{1}{(2^ j\gamma(x))^n}\int_{|x-y|<2^{j+1}\gamma(x)}|f(y)|dy\leq C||f||_{{\rm BMO}^\mathcal{L}(\mathbb{R}^n)},\quad x\in Q_k^*.
\end{eqnarray*}
Also, by using \cite[Lemma 3, (b)]{BHS2}, we obtain
\begin{eqnarray*}
|F_{2,k}(x)|&\leq&V_\rho (R_{\ell,k}^{\mathcal{L},\varepsilon}-R_{\ell,k}^\varepsilon)(f)(x)\\
&=&\sup_{\{\varepsilon _j\}_{j\in \mathbb{N}}\downarrow 0}\left(\sum_{j=1}^\infty \left|\int_{\varepsilon _{j+1}<|x-y|<\varepsilon _j,\,|x-y|<\gamma(x_k)}(R_\ell^\mathcal{L}(x,y)-R_\ell (x-y))f(y)dy\right|^\rho \right)^{1/\rho }\\
&\leq&\int_{|x-y|<\gamma(x_k)}|R_\ell^\mathcal{L}(x,y)-R_\ell(x-y)||f(y)|dy\\
&\leq&C\int_{|x-y|<\gamma(x_k)}\frac{1}{|x-y|^n}\left(\frac{|x-y|}{\gamma(x)}\right)^{2-n/q}|f(y)|dy,\quad x\in Q_k^*.
\end{eqnarray*}

Then, since $\gamma(x)\sim\gamma(x_k)$,
$x\in Q_k^*$, H\"older inequality implies that
\begin{eqnarray*}
|F_{2,k}(x)|&\leq&C\left(\int_{|x-y|<\gamma(x_k)}|x-y|^{(2-n/q-n)r}dy\right)^{1/r}\frac{1}{\gamma(x_k)^{2-n/q}}\left(\int_{|x-y|<\gamma(x_k)}|f(y)|^{r'}dy\right)^{1/r'}\\
&\leq&C\left(\frac{1}{\gamma(x_k)^n}\int_{|x-y|<\gamma
(x_k)}|f(y)|^{r'}dy\right)^{1/r'}\leq C||f||_{{\rm
BMO}^\mathcal{L}(\mathbb{R}^n)},\quad x\in Q_k^*.
\end{eqnarray*}
Here, $1<r<n/(n-2+n/q)$.

Since, for $i=1,2$, $F_{i,k}\in L^\infty (Q_k^*)$ and
$||F_{i,k}||_{L^\infty (Q_k^*)}\leq C||f||_{{\rm
BMO}^\mathcal{L}(\mathbb{R}^n)}$, where $C$ does not depend on $k\in
\mathbb{N}$. In order to see that the properties $(i_k)$ and
$(ii_k)$ hold for $H=V_\rho (R_\ell^{\mathcal{L},\varepsilon})$ it
is sufficient to establish those properties for $H=V_\rho
(R_{\ell,k}^\varepsilon)$.

Fix again $k\in \mathbb{N}$. Then
\begin{eqnarray*}
V_\rho (R_{\ell,k}^\varepsilon)(f)(x)&=&\sup_{\{\varepsilon _j\}_{j\in \mathbb{N}}\downarrow 0}\left(\sum_{j=1}^\infty \left|\int_{\varepsilon _{j+1}<|x-y|<\varepsilon _j,|x-y|< \gamma(x_k)}R_\ell(x-y)f(y)dy\right|^\rho \right)^{1/\rho }\\
&=&\sup_{\{\varepsilon _j\}_{j\in \mathbb{N}}\downarrow 0}\left(\sum_{j=1}^\infty \left|\int_{\varepsilon _{j+1}<|x-y|<\varepsilon _j,|x-y|<\gamma(x_k)}R_\ell(x-y)f(y)\chi _{Q_k^{**}}(y)dy\right|^\rho \right)^{1/\rho }\\
&\leq&V_\rho(R_{\ell,k}^\varepsilon )(f\chi _{Q_k^{**}})(x),\quad
x\in Q_k^*.
\end{eqnarray*}
Hence, according to Theorem \ref{Th2Lp}, we have
\begin{eqnarray}\label{F13}
\frac{1}{|Q_k|}\int_{Q_k}|V_\rho
(R_{\ell,k}^\varepsilon)(f)(x)|dx&\leq&
\left(\frac{1}{|Q_k|}\int_{Q_k}|V_\rho (R_{\ell,k}^\varepsilon)(f\chi _{Q_k^{*}})(x)|^2dx\right)^{1/2}\nonumber\\
&\leq&C\left(\frac{1}{|Q_k|}\int_{Q_k^{*}}|f(x)|^2dx\right)^{1/2}\leq
C||f||_{{\rm BMO}^\mathcal{L}(\mathbb{R}^n)}.
\end{eqnarray}

Let now $x_0\in \mathbb{R}^n$ and $r_0>0$ such that
$B=B(x_0,r_0)\subset Q_k^*$. Then, by using Proposition \ref{Rtrun} we can see
$$
\frac{1}{|B|}\int_{B}|V_\rho (R_{\ell,k}^\varepsilon)(f)(x)-V_\rho
(R_{\ell,k}^\varepsilon)(f_2)(z_1)|dx\leq C||f||_{{\rm
BMO}^\mathcal{L}(\mathbb{R}^n)}.
$$
where $f_2=(f-f_B)\chi _{(B^{**})^c}$ and $z_1\in B$ is such that
$V_\rho (R_{\ell,k}^\varepsilon )(f_2)(z_1)<\infty$. Hence, $V_\rho
(R_{\ell,k} ^\varepsilon )(f)\in {\rm BMO}(Q_k^*)$ and
\begin{equation}\label{F14}
||V_\rho (R_{\ell,k} ^\varepsilon )(f)||_{{\rm BMO}(Q_k^*)}\leq
C||f||_{{\rm BMO}^\mathcal{L}(\mathbb{R}^n)}.
\end{equation}

Note that the constants $C>0$ appearing in (\ref{F13}) and
(\ref{F14}) do not depend on $k\in \mathbb{N}$. Thus the proof of
the desired result is finished.

 $(ii)$ Assume now that $V\in RH_q$ where $n/2<q$. It is sufficient to
consider $n/2<q<n$ (\cite{G}). By proceeding as in the
proof of the previous case, this result will be established when we see
that, for every $k\in \mathbb{N}$, the operators defined by
$$
T_{1,k}(f)(x)=\int_{|x-y|>\gamma(x_k)}|R_\ell^\mathcal{L}(y,x)||f(y)|dy,
$$
and
$$
T_{2,k}(f)(x)=\int_{|x-y|<\gamma(x_k)}|R_\ell^\mathcal{L}(y,x)-R_\ell(y-x)||f(y)|dy,
$$
map ${\rm BMO}^\mathcal{L}(\mathbb{R}^n)$ into $L^\infty (Q_k^*)$, and, for $i=1,2$,
$$
||T_{i,k}(f)||_{L^\infty (Q_k^*)}\leq C||f||_{{\rm BMO}^\mathcal{L}(\mathbb{R}^n)},\quad f\in {\rm BMO}^\mathcal{L}(\mathbb{R}^n),
$$
where $C>0$ does not depend on $k\in \mathbb{N}$.

Let $k\in \mathbb{N}$ and $f\in {\rm BMO}^\mathcal{L}(\mathbb{R}^n)$. According to \cite[p. 538]{Sh}, we have that
\begin{eqnarray*}
|T_{1,k}(f)(x)|&\leq&C\left(\int_{|x-y|>\gamma(x_k)}\frac{1}{|x-y|^n}\frac{1}{(1+|x-y|/\gamma(x))^\alpha}|f(y)|dy\right.\\
&&+\left.\int_{|x-y|>\gamma(x_k)}\frac{1}{|x-y|^{n-1}}\frac{|f(y)|}{(1+|x-y|/\gamma(x))^\alpha}\int_{|y-z|<\frac{|x-y|}{4}}\frac{V(z)}{|z-y|^{n-1}}dzdy\right)\\
&=&C(T _{1,1,k}(f)(x)+T_{1,2,k}(f)(x)),\quad x\in Q_k^*,
\end{eqnarray*}
where $\alpha>0$ will be chosen later large enough.

As it was shown earlier, we have
\begin{equation}\label{F15}
||T_{1,1,k}(f)||_{L^\infty(Q_k^*)}\leq C||f||_{{\rm
BMO}^\mathcal{L}(\mathbb{R}^n)},
\end{equation}
provided that $\alpha\ge 1$.

On the other hand, since $\gamma(x)\sim \gamma(x_k)$ when $x\in Q_k^*$, we can write
\begin{eqnarray*}
|T_{1,2,k}(f)(x)|&\leq&C\sum_{j=0}^\infty \frac{1}{2^{j\alpha }(2^j\gamma(x_k))^{n-1}}\int_{2^j\gamma(x_k)<|x-y|\leq 2^{j+1}\gamma(x_k)}|f(y)|\int_{B(y,\frac{|x-y|}{4})}\frac{V(z)}{|z-y|^{n-1}}dzdy\\
&\leq&C\sum_{j=0}^\infty \frac{1}{2^{j\alpha }(2^j\gamma(x_k))^{n-1}}\left(\int_{|x-y|\leq 2^{j+1}\gamma(x_k)}|f(y)|^{p_0'}dy\right)^{1/p_0'}\\
&&\times\left(\int_{\mathbb{R}^n}\left|\int_{|x-z|<2^{j+2}\gamma(x_k)}\frac{V(z)}{|z-y|^{n-1}}dz\right|^{p_0}dy\right)^{1/p_0},
\quad x\in Q_k^*,
\end{eqnarray*}
where $\frac{1}{p_0}=\frac{1}{q}-\frac{1}{n}$. Then, the
$L^p$-boundedness properties of the fractional integrals
(\cite[p. 354]{Stein1}) lead us to
$$
|T_{1,2,k}(f)(x)|\leq C||f||_{{\rm
BMO}^\mathcal{L}(\mathbb{R}^n)}\sum_{j=0}^\infty \frac{1}{2^{j\alpha
}(2^j\gamma(x_k))^{n-1-n/q_0'}}\left(\int_{|x-z|<2^{j+2}\gamma(x_k)}|V(z)|^qdz\right)^{1/q},
\;x\in Q_k^*.
$$
By using the properties of $V$ and $\gamma$ (\cite[Lemma 1]{BHS2}) we
obtain, for a certain $\mu >0$,
$$
\left(\int_{|x-z|<2^{j+2}\gamma(x_k)}|V(z)|^qdz\right)^{1/q}\leq
C(2^j\gamma(x_k))^{-n/q'}2^{j\mu }\gamma(x_k)^{n-2},\quad x\in Q_k^*.
$$
It follows that, by choosing $\alpha >0$ large enough,
\begin{eqnarray}\label{F16}
|T_{1,2,k}(f)(x)|&\leq&C||f||_{{\rm BMO}^\mathcal{L}(\mathbb{R}^n)}\sum_{j=0}^\infty \frac{1}{2^{j(\alpha +n/q_0-1+n/q'-\mu)}}\nonumber\\
&\leq&C||f||_{{\rm BMO}^\mathcal{L}(\mathbb{R}^n)},\quad x\in Q_k^*.
\end{eqnarray}
We conclude from (\ref{F15}) and (\ref{F16}) that
$$
||T_{1,k}(f)||_{L^\infty (Q_k^*)}\leq C||f||_{{\rm
BMO}^\mathcal{L}(\mathbb{R}^n)},
$$
where $C>0$ does not depend on $k\in \mathbb{N}$.

According to \cite[(5.9)]{Sh} we get
\begin{eqnarray*}
|T_{2,k}(f)(x)|&\leq&C\left(\int_{|x-y|<\gamma(x_k)}\frac{1}{|x-y|^n}\left(\frac{|x-y|}{\gamma(x)}\right)^{2-n/q}|f(y)|dy\right.\\
&+&\left.\int_{|x-y|<\gamma
(x_k)}\frac{1}{|x-y|^{n-1}}\int_{|y-z|<\frac{|x-y|}{4}}\frac{V(z)}{|z-y|^{n-1}}dz|f(y)|dy\right)\\
&=&C(T_{2,1,k}(f)(x)+T_{2,2,k}(f)(x)),\,\,\,x\in Q_k^*.
\end{eqnarray*}
By proceeding as in the first part of this proof we have that
\begin{equation}\label{F17}
||T_{2,1,k}(f)||_{L^\infty (Q_k^*)}\leq C||f||_{{\rm
BMO}^\mathcal{L}(\mathbb{R}^n)}.
\end{equation}
We can also write
\begin{eqnarray*}
|T_{2,2,k}(f)(x)|&\leq&C\sum_{j=0}^\infty \int_{2^{-j-1}\gamma(x_k)\leq |x-y|<2^{-j}\gamma(x_k)}\frac{|f(y)|}{(2^{-j}\gamma(x_k))^{n-1}}\int_{|x-z|<2^{-j+1}\gamma(x_k)}\frac{V(z)}{|y-z|^{n-1}}dzdy\\
&\leq&C\sum_{j=0}^\infty \frac{1}{(2^{-j}\gamma(x_k))^{n-1}}\left(\int_{|x-y|<2^{-j}\gamma(x_k)}|f(y)|^{p_0'}dy\right)^{1/p_0'}\\
&&\times\left(\int_{\mathbb{R}^n}\left(\int_{|x-z|<2^{-j+1}\gamma(x_k)}\frac{V(z)}{|y-z|^{n-1}}dz\right)^{p_0}dy\right)^{1/p_0}\\
&\leq&C\sum_{j=0}^\infty \frac{1}{(2^{-j}\gamma(x_k))^{n-1-n/p_0'}}\left(\int_{|x-z|<2^{-j+1}\gamma(x_k)}V(z)^qdz\right)^{1/q}\\
&\times&\left(\frac{1}{(2^{-j}\gamma(x_k))^n}\int_{|x-y|<2^{-j}\gamma(x_k)}|f(y)|^{p_0'}dy\right)^{1/p_0'},\quad x\in Q_k^*,
\end{eqnarray*}
where $\frac{1}{p_0}=\frac{1}{q}-\frac{1}{n}$.

Since $V\in RH_q$ and $\gamma(x)\sim \gamma(x_k)$, when $x\in
Q_k^*$, we have (\cite[Lemma 1]{BHS2})
$$
\left(\int_{B(x,2^{-j+1}\gamma
(x_k))}V(z)^qdz\right)^{1/q}\leq C\gamma(x_k)^{n/q-2},\quad x\in
Q_k^*.
$$
Moreover, an argument like the one used to show \cite[Lemma 2]{DGMTZ}
allows us to get
$$
\left(\frac{1}{(2^{-j}\gamma(x_k))^n}\int_{|x-y|<2^{-j}\gamma(x_k)}|f(y)|^{p_0'}dy\right)^{1/p_0'}\leq Cj||f||_{{\rm BMO}^\mathcal{L}(\mathbb{R}^n)}.
$$
Then,
\begin{eqnarray}\label{F18}
|T_{2,2,k}(f)(x)|&\leq&C\sum_{j=0}^\infty \frac{j}{(2^{-j}\gamma(x_k))^{n/p_0-1}}\gamma(x_k)^{n/q-2}||f||_{{\rm BMO}^\mathcal{L}(\mathbb{R}^n)}\nonumber\\
&\leq&C||f||_{{\rm BMO}^\mathcal{L}(\mathbb{R}^n)},\quad x\in Q_k^*.
\end{eqnarray}
Note that $\frac{n}{p_0}-1=\frac{n}{q}-2<0$.

By combining (\ref{F17}) and (\ref{F18}) we conclude that
$$
||T_{2,k}(f)||_{L^\infty (Q_k^*)}\leq C||f||_{{\rm
BMO}^\mathcal{L}(\mathbb{R}^n)},
$$
where $C>0$ does not depend on $k\in \mathbb{N}$.

Thus the proof is finished.

  \end{document}